%% file: bad.tex
\documentclass{amsart}

\title{Three point covers with bad reduction}

\author{Stefan Wewers}

\address{Mathematisches Institut\\
         Beringstr. 1\\
         53115 Bonn, Germany}

\email{wewers@math.uni-bonn.de}

\date{}

\subjclass[2000]{Primary 14H30, 11G20}

\keywords{Three point cover, stable reduction, field of moduli}

\usepackage{amssymb}
\usepackage{latexsym}
\usepackage{amscd}

\theoremstyle{plain}

        \newtheorem{thm}{Theorem}[section]
        \newtheorem{cor}[thm]{Corollary}
        \newtheorem{lem}[thm]{Lemma}
        \newtheorem{prop}[thm]{Proposition}
        \newtheorem{introthm}{Theorem}

\theoremstyle{definition}

        \newtheorem{defn}[thm]{Definition}
        \newtheorem{ass}[thm]{Assumption}
        
        \newtheorem{exa}[thm]{Example}

\theoremstyle{remark}

        \newtheorem{rem}[thm]{Remark}

\newcommand{\ZZ}{\mathbb{Z}}
\newcommand{\QQ}{\mathbb{Q}}
\newcommand{\CC}{\mathbb{C}}
\newcommand{\PP}{\mathbb{P}}

\newcommand{\FF}{\mathbb{F}}

\newcommand{\C}{\mathfrak{C}}

\newcommand{\OO}{\mathcal{O}}
\newcommand{\B}{\mathcal{B}}

\newcommand{\D}{\mathcal{D}}
\newcommand{\U}{\mathcal{U}}
\newcommand{\V}{\mathcal{V}}
\newcommand{\E}{\mathcal{E}}

\newcommand{\Od}{\hat{\OO}}
\newcommand{\Z}{\mathcal{Z}}
\newcommand{\G}{\mathcal{G}}

\newcommand{\T}{\mathcal T}

\newcommand{\X}{\mathcal{X}}
\newcommand{\Y}{\mathcal{Y}}

\newcommand{\Spec}{\mathop{\rm Spec}}

\newcommand{\EXt}{\mathop{\mathcal{E}xt}\nolimits} 
\newcommand{\Ind}{\mathop{\rm Ind}\nolimits} 
 
\newcommand{\Gal}{\mathop{\rm Gal}} 
\newcommand{\Aut}{\mathop{\rm Aut}\nolimits} 
\newcommand{\ord}{\mathop{\rm ord}\nolimits} 
\newcommand{\val}{\mathop{\rm val}\nolimits}

\newcommand{\Ker}{\mathop{\rm Ker}} 
\newcommand{\Def}{\mathop{\rm Def}\,}

\newcommand{\inn}{^{\rm\scriptscriptstyle in}}
\newcommand{\ab}{^{\rm\scriptscriptstyle ab}}
\newcommand{\stab}{^{\rm\scriptscriptstyle st}}
\newcommand{\aux}{^{\rm\scriptscriptstyle aux}}
\newcommand{\dagg}{^{\scriptscriptstyle\dagger}}
\newcommand{\prim}{_{\rm\scriptscriptstyle prim}}
\newcommand{\new}{_{\rm\scriptscriptstyle new}}

\newcommand{\tame}{_{\rm\scriptscriptstyle tame}}
\newcommand{\wild}{_{\rm\scriptscriptstyle wild}}
\newcommand{\branch}{_{\rm\scriptscriptstyle ram}}
\newcommand{\intern}{_{\rm\scriptscriptstyle int}}
\newcommand{\diff}{{\rm d}}
\newcommand{\infb}{\bar{\infty}}

\newcommand{\fb}{\bar{f}}
\newcommand{\gb}{\bar{g}}

\newcommand{\Rt}{\tilde{R}}
\newcommand{\Xb}{\bar{X}}
\newcommand{\Yb}{\bar{Y}}
\newcommand{\Zb}{\bar{Z}}
\newcommand{\Wb}{\bar{W}}
\newcommand{\Xd}{\hat{X}}
\newcommand{\Yd}{\hat{Y}}
\newcommand{\Zd}{\hat{Z}}

\newcommand{\xb}{\bar{x}}

\newcommand{\ub}{\bar{u}}

\newcommand{\Ub}{\bar{U}}
\newcommand{\Vb}{\bar{V}}
\newcommand{\wb}{\bar{w}}

\newcommand{\psib}{\bar{\psi}}

\newcommand{\vphi}{\varphi}
\newcommand{\vphib}{\bar{\vphi}}

\newcommand{\kappab}{\bar{\kappa}}

\newcommand{\inj}{\hookrightarrow}
\newcommand{\To}{\;\longrightarrow\;}
\newcommand{\iso}{\stackrel{\sim}{\to}}
\newcommand{\liso}{\;\stackrel{\sim}{\longrightarrow}\;}

\newcommand{\lpfeil}[1]{\stackrel{#1}{\To}}

\newcommand{\op}[2]{\sideset{^{#1}}{}{\mathop{#2}}}

\newcommand{\p}{\mathfrak{p}}

\newcommand{\bmu}{\boldsymbol{\mu}}
\newcommand{\balpha}{\boldsymbol{\alpha}}

\newcommand{\Kb}{\bar{K}}

\newcommand{\gen}[1]{\mathopen\langle#1\mathclose\rangle}
\newcommand{\vf}[1]{\frac{\diff}{\diff #1}}

\begin{document}

\begin{abstract}
  We study Galois covers of the projective line branched at three
  points with bad reduction to characteristic $p$, under the condition
  that $p$ strictly divides the order of the Galois group.  As an
  application of our results, we prove that the field of moduli of
  such a cover is at most tamely ramified at $p$.
\end{abstract}

\maketitle

\section*{Introduction}

\subsection*{Ramified primes in the field of moduli}

Let $f:Y\to X=\PP^1_{\CC}$ be a three point cover, i.e.\ a finite
cover of the Riemann sphere branched at $0$, $1$ and $\infty$. The
`obvious direction' of Belyi's Theorem states that $f$ can be defined
over a number field. Therefore, the absolute Galois group of $\QQ$
acts on the set of isomorphism classes of all three point covers, and
we can associate to $f$ the number field $K$ such that
$\Gal(\bar{\QQ}/K)$ is the stabilizer of the isomorphism class of $f$.
The number field $K$ is called the field of moduli of $f$. Very little
is known about the correspondence between three point covers and their
associated field of moduli. One of the most general things that we
know about $K$ is given by the following theorem of Beckmann
\cite{Beckmann89}.  Let $G$ be the monodromy group of $f$. Then the
extension $K/\QQ$ is unramified outside the set of primes diving the
order of $G$. This result is related to the fact that $f$ has good
reduction at each prime ideal of $K$ dividing such primes $p$.

A recent theorem of Raynaud \cite{Raynaud98} gives a partial converse to
Beckmann's Theorem. Suppose that the prime $p$ {\em strictly} divides the
order of $G$, i.e.\ that $p^2$ does not divide $|G|$. Let $\p$ be a prime
ideal of $K$ dividing $p$ and suppose that the ramification index $e(\p/p)$ is
strictly smaller than the number of conjugacy classes in $G$ of elements of
order $p$. Then $f$ has potentially good reduction at $\p$.

The proof of this theorem depends on an analysis of the {\em stable reduction}
of $f:Y\to X$ at $\p$. If $f$ does not have potentially good reduction at $\p$
then this analysis yields a certain lower bound for the ramification index
$e(\p/p)$. The theorem follows from this bound. We remark that this sort of
result is not restricted to three point covers, but extends, under certain
conditions, to covers of more general curves $X$. The essential condition here
is that $p$ strictly divides the order of $G$.

In \cite{special} the author has continued Raynaud's study of the
stable reduction, and the present paper is a further continuation of
this work. Our results yield an {\em upper} bound for the ramification
index of $p$ in the field of moduli of certain three point covers. For
instance, we prove the following theorem.

\begin{introthm} \label{introthm1}
  Let $f:Y\to\PP^1$ be a three point cover, with monodromy group $G$
  and field of moduli $K$.  Let $p$ be a prime which {\em strictly}
  divides the order of $G$. Then $p$ is at most tamely ramified in the
  extension $K/\QQ$.
\end{introthm}

In fact we obtain a formula for the ramification index
$e(\p/p)$ in terms of the stable reduction of $f:Y\to X$ at $\p$. This
formula shows that $e(\p/p)$ is prime to $p$ and gives an upper bound.

In contrast to the theorem of Raynaud mentioned above, our results
are very particular to three point covers. They depend in an essential
way on the fact that three point covers are `rigid' objects, i.e.\ do
not admit any nontrivial deformation.

\subsection*{Stable reduction}

Let $R$ be a complete discrete valuation ring, with residue field $k$
of characteristic $p>0$ and fraction field $K$ of characteristic zero.
Let $f_K:Y_K\to X_K:=\PP^1_K$ be a three point cover defined over $K$.
We assume that $f_K$ is Galois, with Galois group $G$. In order to
prove results such as Theorem \ref{introthm1}, this is no restriction.
After replacing the field $K$ by a finite extension, we may assume
that $Y_K$ has a semistable model $Y_R$ over $R$ such that the
ramification points of $f_K$ specialize to pairwise distinct smooth
points on the special fiber. Among all semistable models with this
property, we let $Y_R$ be the minimal one. The action of $G$ extends
to $Y_R$; let $X_R:=Y_R/G$ denote the quotient.  This is a
semistable model of $\PP^1_K$ over $R$. Therefore, the special fiber
$\Xb:=X_R\otimes_R k$ is a tree of curves of genus zero.  The map
$f_R:Y_R\to X_R$ is called the {\em stable model} of the Galois cover
$f_K$. Its restriction to the special fiber $\fb:\Yb\to\Xb$ is called
the {\em stable reduction} of $f_K$.

We say that the cover $f_K:Y_K\to X_K$ has good reduction if the
curves $X_R$ and $Y_R$ are smooth over $R$. If this is the case, then
the stable reduction $\fb:\Yb\to\Xb=\PP^1_k$ is a three point Galois
cover in characteristic $p$. If $f_K$ does not have good reduction we
say that it has {\em bad reduction}.  Our first main result is a
structure theorem for the stable reduction $\fb:\Yb\to\Xb$ in the case
of bad reduction.

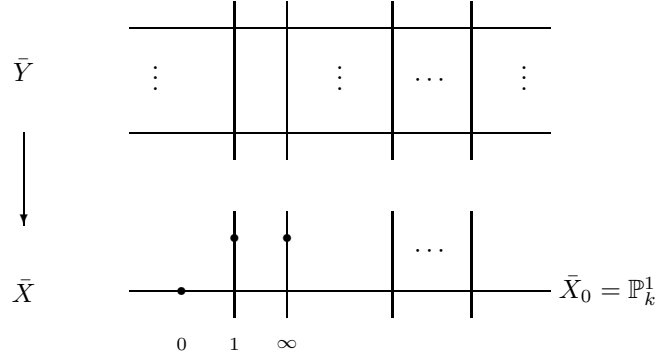
\begin{figure}
\begin{center}
\unitlength3.5mm
\input{stablepic.tex}
\caption{\label{stablepic} The stable reduction of a three point cover}
\end{center}
\end{figure}

\begin{introthm} \label{introthm2}
  Suppose that $p$ strictly divides the order of $G$ and that the
  cover $f_K:Y_K\to X_K$ has bad reduction. Then the following holds
  (compare with Figure \ref{stablepic}).
  \begin{enumerate}
  \item The curve $\Xb$ consists of one central component $\Xb_0$,
    canonically isomorphic to $\PP^1_k$, and $r\geq 1$ tails
    $\Xb_1,\ldots,\Xb_r$, each of which intersects $\Xb_0$ in one
    point.
  \item
    Let $\Yb_0$ be a irreducible component of $\Yb$ lying above the
    central component $\Xb_0$. Then the map $\Yb_0\to\Xb_0$ induced
    from $\fb$ is the composition of a purely inseparable map
    $\Yb_0\to\Zb_0$ of degree $p$ and a Galois cover $\Zb_0\to\Xb_0$.
  \item For $j=1,\ldots,r$, let $\Yb_j$ be an irreducible component of
    $\Yb$ lying over the tail $\Xb_j$. Then the map $\Yb_j\to\Xb_j$
    induced by $\fb$ is a Galois cover, ramified at most at two
    points. More precisely, $\Yb_j\to\Xb_j$ is wildly ramified at the
    point where $\Xb_j$ intersects $\Xb_0$ and is tamely ramified at
    any point which is the specialization of one of the $K$-points
    $0$, $1$ or $\infty$.
  \end{enumerate}
\end{introthm}

Part (1) is the essential statement and depends on the assumption that
$f_K$ is a three point cover. Part (2) and (3) follow from (1), using
the results of \cite{Raynaud98} and the assumption that $p$ strictly
divides the order of $G$.

Under the additional assumption that all the ramification indices of
$f_K$ are prime to $p$, Theorem \ref{introthm2} (1) can be deduced
from the results of \cite{special}, via Raynaud's construction of the
{\em auxiliary cover} (see the introduction of \cite{special}). The
proof of the general case given in the present paper is based on a
generalization of the methods of \cite{special}.  This generalization
avoids the assumption on the ramification indices and the use of the
auxiliary cover.

\subsection*{Lifting of special $G$-maps}

The stable reduction of a three point Galois cover $f_K:Y_K\to
X_K:=\PP^1_K$ is, by definition, a finite map $\fb:\Yb\to\Xb$ between
semistable curves over the residue field $k$, together with an action
of a finite group $G$ on $\Yb$ which commutes with $\fb$. In the case
of bad reduction, the curves $\Xb$ and $\Yb$ are both singular, and
the map $\fb$ is inseparable over some of the irreducible components
of $\Xb$. This suggests the following question. Given a map
$\fb:\Yb\to\Xb$ of the sort we have just described, does it occur as
the stable reduction of a three point Galois cover $f_K:Y_K\to X_K$ in
characteristic zero? If this is the case then we say that $f_K:Y_K\to
X_K$ is a {\em lift} of $\fb:\Yb\to\Xb$.

Theorem \ref{introthm2} gives a list of necessary conditions on $\fb$ for the
existence of a lift (under the condition that $p$ strictly divides the order
of the group $G$). There is one extra condition on $\fb$ which we have omitted
from the statement of Theorem \ref{introthm2}, but which plays a central role
in this paper. This additional condition is the existence of a certain
differential form $\omega_0$ on the curve $\Zb_0$ which is `compatible' with
$\fb$ (see Theorem \ref{introthm2} (2) for the definition of $\Zb_0$). For
instance, $\omega_0$ is {\em logarithmic}, i.e.\ of the form $\diff u/u$, and
it is an eigenvector under the action of the Galois group of the cover
$\Zb_0\to\Xb_0$. Moreover, the divisor of zeroes and poles of $\omega_0$ is
related to and determined by the map $\fb$ and the action of $G$ on $\Yb$. The
existence and the basic properties of $\omega_0$ follow from work of
Green--Matignon \cite{GreenMatignon99} and Henrio \cite{YannickArbres}. In
some sense, $\omega_0$ carries infinitesimal information about the action of
$G$ on $\Yb$ which `got lost' by reduction to characteristic $p$. Following
\cite{special}, we call $(\Zb_0,\omega_0)$ a {\em special deformation datum}.
(Actually, in \cite{special} we assume that $\omega_0$ has no poles.  This is
true if and only if all the ramification indices of the three point cover
$f_K:Y_K\to X_K$ are prime to $p$.)

In the present paper we define the notion of a {\em special $G$-map}.
This is a map $\fb:\Yb\to\Xb$ between semistable curves in
characteristic $p$ together with an action of a finite group $G$ on
$\Yb$ satisfying certain conditions. Essentially, these conditions say
that $\fb$ satisfies the conclusion of Theorem \ref{introthm2} and
that there exists a special deformation datum $(\Zb_0,\omega)$ which
is `compatibel' with $\fb:\Yb\to\Xb$ and the action of $G$ on $\Yb$.
Our second main result can be stated as follows.

\begin{introthm} \label{introthm3}
  Let $\fb:\Yb\to\Xb$ be a special $G$-map, defined over an
  algebraically closed field $k$ of characteristic $p>0$. Let $K_0$
  denote the fraction field of the ring of Witt vectors over $k$. 
  Then the following holds.
  \begin{enumerate}
  \item There exists a lift of $\fb$, i.e.\ a three point Galois cover
    $f_K:Y_K\to X_K$ whose stable reduction is isomorphic to
    $\fb:\Yb\to\Xb$.
  \item
    Each lift of $\fb:\Yb\to\Xb$ can be defined over a finite extension
    $K/K_0$ which is at most tamely ramified. 
  \end{enumerate}
\end{introthm} 

We actually prove more, namely we determine the set of isomorphism
classes of all lifts of $\fb$, together with the action of
$\Gal(\Kb_0/K_0)$ on this set, in terms of certain numerical invariants
attached to $\fb$. From this more precise result one can deduce an
upper bound for the degree of the minimal extension $K/K_0$ over which
a given lift of $\fb$ can be defined. In any case, Theorem \ref{introthm1}
follows from Theorem \ref{introthm2} and Theorem \ref{introthm3} by a
straightforward argument.

Part (1) of Theorem \ref{introthm3}, i.e.\ the mere existence of a
lift, follows already from the results of \cite{special} (assuming
that $\omega_0$ has no poles).  Part (2) is more difficult and needs
new methods. The key step in the proof of (2) is the study of the
deformation theory of a certain curve with group scheme action
associated to the special deformation datum $(\Zb_0,\omega_0)$. A
detailed exposition of this deformation theory can be found in
\cite{cotang}.

\vspace{1ex} This paper is divided into four sections. The first
section contains the general theory of stable reduction of Galois
covers, the monodromy action and deformation data. In \S
\ref{threepoints} we apply this theory to the case of three point
cover and prove our first main result, which essentially corresponds
to Theorem \ref{introthm1} above. We also define the notion of a
special $G$-map. In \S \ref{aux} we apply the deformation theory of
\cite{cotang} in order to lift a special deformation datum to an
auxiliary cover with certain good properties. This is preparatory work
for \S \ref{main}, where we prove the lifting result corresponding to
Theorem \ref{introthm3}. At the end of \S \ref{main} we draw some
conclusions from this result, concerning the field of moduli, and we
give some examples.

\subsection*{Acknowledgments}

This work grew out of discussions with Michel Raynaud on my previous paper
\cite{special}. I would like to thank him heartily for his critical remarks
and encouragement. The final version of this paper was written while I was a
guest of the Max--Planck--Institut in Bonn.





\section{The stable reduction of Galois covers} \label{stablered}

Let $K$ be a field of characteristic zero equipped with a complete
discrete valuation whose residue field $k$ is algebraically closed and
has characteristic $p>0$. Let $f:Y\to X$ be a Galois cover of smooth
projective curves over $K$. In this section we define and study the
{\em stable reduction} of $f$. We make two crucial assumptions on the
reduction of $f$. First, we assume that the marked curve $(X,S)$
(where $S$ is the branch locus of $f$) has good reduction. Second, we
assume that $f$ has {\em mild reduction}, a property which is
automatically verified if $p$ strictly divides the order of the Galois
group of $f$.

We start by recalling some results of Raynaud \cite{Raynaud98} on the
structure of $\fb:\Yb\to\Xb$. Then we attach to each irreducible
component of $\Xb$ where the stable reduction $\fb:\Yb\to\Xb$ is
inseparable a so-called {\em deformation datum}. This datum encodes
infinitesimal information on the action of the Galois group of $f:Y\to
X$. The deformation data attached to the individual components of
$\Xb$ are related by certain compatibility conditions. The existence
of a set of compatible deformation data imposes strong restrictions on
the map $\fb:\Yb\to\Xb$. In this section, we content ourself with
working out those restrictions which are of a purely combinatorial
nature. For instance, we obtain a new proof of Raynaud's {\em
  vanishing cycle formula} \cite{Raynaud98}.

Most of the results of this section are already present in
\cite{special}, in a special case.

\subsection{The stable model and the monodromy action} \label{stable}

\subsubsection{} \label{stable1} 

The following notation is fixed throughout \S \ref{stablered} and
\S \ref{threepoints}.  Let $R_0$ be a complete discrete valuation
ring, with quotient field $K_0$ of characteristic $0$ and residue
field $k$ of characteristic $p>0$.  We assume that $k$ is
algebraically closed. We fix an algebraic closure $\Kb$ of $K_0$.

Let $(X_0,S_0)$ be a {\em smooth} stably marked curve over $R_0$.
Since $R_0$ is strictly henselian this implies that
$S_0=\{\,x_{R_0,j}\mid j\in B_0\}$, where $x_{R_0,j}:\Spec R_0\to X_0$
are pairwise distinct sections, indexed by a finite set $B_0$. If
$K/K_0$ is a field extension then we write $X_K$ and $S_K$ instead of
$X_0\otimes_{R_0}K$ and $S_0\otimes_{R_0}K$. If $K=\Kb$ then we simply
write $X$ and $S$ instead of $X_{\Kb}$ or $S_{\Kb}$. The elements of
$S$ are denoted by $x_j$, $j\in B_0$. In \S \ref{threepoints} we will
consider the case $X_0=\PP^1_{R_0}$ and $S_0=\{0,1,\infty\}$.

Let $G$ be a finite group and $f:Y\to X$ a $G$-cover of $X$ with
branch locus $S$. We consider the set of ramification points of the
cover $f$ as a marking on $Y$; this makes $Y$ a stably marked curve.
Let $K$ be a finite extension of $K_0$ and $R$ the ring of integers of
$K$. A {\em model} of $f$ over $K$ is a $G$-cover $f_K:Y_K\to X_K$
such that $f=f_K\otimes_K\Kb$.  A {\em good model} of $f$ over $K$ is
a model which extends to a finite and tame cover $f_R:Y_R\to
X_0\otimes_{R_0}R$, ramified along $S_0\otimes_{R_0}R$ and \'etale
everywhere else.  We say that $f$ has {\em good reduction} if there
exists a good model over some extension $K/K_0$.  If $f$ does not have
good reduction, we say that it has {\em bad reduction}.

The following is well known, see \cite{Beckmann89}. If $f$ has good
reduction then there exists a good model $f_{K_0}:Y_{K_0}\to X_{K_0}$
of $f$ over $K_0$, unique up to isomorphism. Moreover, if $p$ does not
divide the order of the Galois group $G$, then $f:Y\to X$ has good
reduction.

\subsubsection{} \label{stable2}

From now on, we assume that $f:Y\to X$ has bad reduction. Note that
this implies that $p$ divides the order of $G$. If $K/K_0$ is a
sufficiently large finite extension, then there exists a model
$f_K:Y_K\to X_K$ of $f$ over $K$ such that $Y_K$ extends to a stably
marked curve $Y_R$ over $R$. In particular, $Y_R$ is a semistable
model of $Y_K$ over $R$ such that the ramification points of the cover
$f_K$ specialize to pairwise disjoint, smooth points on the special
fiber $\Yb:=Y_R\otimes k$.

Since the stably marked model is unique, the action of the group $G$
on $Y_K$ extends to $Y_R$. Let $X_R:=Y_R/G$ be the quotient; then
$X_R$ is a semistable curve over $R$ with generic fiber $X_K$, see
\cite{Raynaudfest}, Appendix. Moreover, the branch points $x_j$ of
$f_K$ specialize to pairwise distinct smooth points $\xb_j$ on the
special fiber $\Xb:=X_R\otimes k$. We consider $\Xb$ as a marked
semistable curve (note, however, that $\Xb$ is {\em not} stably
marked). By \cite{Knudsen83}, there exists a canonical contraction
morphism $q_R:X_R\to X_0\otimes_{R_0}R$. Let $\Xb$ denote the special
fiber of $X_R$. There is a unique irreducible component of $\Xb$ on
which the induced map $\bar{q}:\Xb\to\Xb_0:=X_0\otimes_{R_0}k$ is an
isomorphism. We may and will identify this component with $\Xb_0$. All
other components of $\Xb$ are of genus $0$ and are contracted by
$\bar{q}$ to a closed point of $\Xb_0$.

\begin{defn} \label{stabledef}
  The natural map $f_R:Y_R\to X_R$ is called the {\em stable model} of
  the Galois cover $f_K$. The induced map $\fb:\Yb\to\Xb$ on the
  special fibers is called the {\em stable reduction} of $f_K$. The
  component $\Xb_0$ of $\Xb$ is called the {\em original component}.
\end{defn} 

It is clear that the stable model is `stable' under extension of the
field $K$. In particular, the stable reduction does not depend on the
choice of $K$. We remark that the stable model defined above differs
slightly from the model used in \cite{Raynaud98}. This difference is
not essential, but certain results of \cite{Raynaud98} have to be
reformulated. This is done in Lemma \ref{stablem} below.

\subsubsection{} \label{monodromy}

For an element $\sigma\in\Gal(\Kb/K_0)$ we denote by
$\op{\sigma}{f}:\op{\sigma}{Y}\to\op{\sigma}{X}=X$ the conjugate
$G$-cover. Let $\Gamma\inn$ be the subgroup of $\Gal(\Kb/K_0)$
consisting of all elements $\sigma$ such that $f\cong\op{\sigma}{f}$,
as $G$-covers. Similarly, let $\Gamma\ab$ be the subgroup of
$\Gal(\Kb/K_0)$ consisting of elements $\sigma$ such that
$\op{\sigma}{f}\cong f$ as mere covers. Set
\[
    K\inn \;:=\; \Kb^{\Gamma\inn}, \qquad K\ab \;:=\; \Kb^{\Gamma\ab}.
\]
The field $K\inn$ (resp.\ $K\ab$) is the {\em field of moduli} of the
$G$-cover $f$ (resp.\ of the mere cover $f$), relative to the
extension $\Kb/K_0$. See \cite{DebDou1} for more details.

Let $K/K_0$ be a Galois extension such that $f$ has stable reduction
over $K$, and let $f_R:Y_R\to X_R$ be the stable model. Also, let
$\sigma\in\Gamma\ab$. By definition of $\Gamma\ab$, there exists a
$\sigma$-linear automorphism $\kappa_{K,\sigma}:Y_K\iso Y_K$ which
normalizes the action of $G$ and induces the canonical $\sigma$-linear
automorphism $X_K\iso X_K$. The automorphism $\kappa_{K,\sigma}$ is
unique up to composition with an element of $G$. It extends to a
$\sigma$-linear automorphism $\kappa_{R,\sigma}:Y_R\iso Y_R$. Let
$\kappab_\sigma:\Yb\iso\Yb$ be the restriction of $\kappa_{R,\sigma}$
to the special fiber. This is a $k$-linear automorphism which
normalizes the action of $G$. It is uniquely determined by
$\sigma\in\Gamma\ab$, up to composition with an element of $G$. We
obtain a continuous homomorphism
\[
     \kappab\ab:\,\Gamma\ab \To \Aut_k(\Yb)/G,
\]
which we shall call the (absolute) {\em monodromy action}. Note that
$\kappab_\sigma$ induces an automorphism of $\Xb$ which is the
identity on the original component. It follows from \cite{DebDou1}
that the set of $K\ab$-models of the mere cover $f$ is in bijection
with the set of homomorphisms $\Gamma\ab\to\Aut_k(\Yb)$ which lift
$\kappab\ab$. If a model $f_{K\ab}$ of $f$ is specified, we will
sometimes consider $\kappab\ab$ as a homomorphism
$\Gamma\ab\to\Aut_k(\Yb)$.

Let $\Aut_G(\Yb)$ denote the group of $k$-linear automorphisms of
$\Yb$ which commute with the $G$-action, and let $C_G$ denote the
center of $G$. The restriction of $\kappab\ab$ to $\Gamma\inn$ yields
a homomorphism
\[
     \kappab\inn:\,\Gamma\inn \To \Aut_G(\Yb)/C_G,
\]
called the inner monodromy action. Similar to what we stated above,
$K\inn$-models of the $G$-cover $f$ correspond one-to-one to 
lifts of $\kappab\inn$ to a homomorphism $\Gamma\inn\to\Aut_G(\Yb)$.

Let $\Gamma\stab$ denote the subgroup of $\Gamma\ab$ of elements
$\sigma$ such that $\kappab\ab_\sigma\in G$, and set
$K\stab:=\Kb^{\Gamma\stab}$. One shows that $K\stab/K_0$ is the
smallest extension over which $f$ has stable reduction. Compare with
\cite{Raynaud98}, \S 2.2.

\subsection{Mild reduction}

\subsubsection{} \label{mild1}

Let $f:Y\to X$ be a $G$-cover as in \S \ref{stable}. We assume that
$f$ has bad reduction and denote by $\fb:\Yb\to\Xb$ its stable
reduction. If $\eta$ is a point on $\Yb$ (closed or generic), we
denote by $G_\eta$ (resp.\ $I_\eta$) the corresponding decomposition
(resp.\ inertia) group. By definition, $G_\eta$ is the set of elements
of $G$ which stabilize $\eta$, and $I_\eta$ is the normal subgroup of
$G_\eta$ consisting of elements which act trivially on the residue
field $k(\eta)$.

\begin{defn} \label{milddef}
  We say that $f:Y\to\PP^1$ has {\em mild reduction} it has bad
  reduction and if the following two conditions hold for all points
  $\eta\in\Yb$. Firstly, the $p$-part of $I_\eta$ has at most order
  $p$. Secondly, if $I_\eta\not=1$ then the order of $G_\eta/I_\eta$
  is prime to $p$.
\end{defn}

If $f:Y\to X$ has bad reduction and $p$ strictly divides the order of
the Galois group $G$, then the reduction of $f$ is automatically mild.
From now on, let us assume that $f:Y\to X$ has mild reduction.

\subsubsection{} \label{mild2}

Let us denote by $(\Xb_v)$ the family of irreducible components of the
curve $\Xb$. 

\begin{defn} \label{taildef}
  A component $\Xb_v$ of $\Xb$ which is not the original
  component and which is connected to the rest of $\Xb$ in a single
  point is called a {\em tail}. The point where a tail $\Xb_v$ meets
  the rest of $\Xb$ is denoted by $\infty_v$. If a tail $\Xb_v$
  contains one of the points $\xb_j$, $j\in B_0$, then $\Xb_v$ is
  called {\em primitive}.  Otherwise, the tail $\Xb_v$ is called {\em
    new}. A component $\Xb_v$ which is not a tail is called an {\em
  interior component}. 
\end{defn}

For each component $\Xb_v$ of $\Xb$, we choose an irreducible
component $\Yb_v$ of $\Yb$ above $\Xb_v$. We write $G_v$ for the
decomposition group of $\Yb_v$ and $I_v$ for the inertia group of
$\Yb_v$.

\begin{lem}[Raynaud] \label{stablem}
  Suppose that $f:Y\to\PP^1$ has mild reduction. Then the
  following holds. 
  \begin{enumerate}
  \item If $\Xb_v$ is a tail, then the inertia group $I_v$ is trivial.
    In other words, the map $\Yb_v\to\Xb_v$ is a Galois cover, with
    Galois group $G_v$. If $\Xb_v$ is a primitive tail and contains
    the point $\xb_j$, $j\in B_0$, then the cover
    $\Yb_v\to\Xb_v$ is wildly ramified at $\infty_v$,
    tamely ramified at $\xb_j$ and \'etale everywhere else. If $\Xb_v$
    is a new tail, then $\Yb_v\to\Xb_v$ is wildly ramified at
    $\infty_v$ and \'etale everywhere else.
  \item If $\Xb_v$ is not a tail, then the inertia group $I_v$ is
    cyclic of order $p$. The cover $\Yb_v\to\Xb_v$ factors through the
    $k$-linear Frobenius map $F:\Yb_v\to\Zb_v:=\Yb_v^{(p)}$.  The
    resulting map $\Zb_v\to\Xb_v$ is a Galois cover, with Galois group
    $H_v:=G_v/I_v$. It is \'etale outside the set of points which are
    either singular points of $\Xb$ or belong to the set $\{\xb_j\mid
    j\in B_0\}$.
  \item For $j\in B_0$, let $n_j$ be the ramification index of $f:Y\to
    X$ in $x_j$. Suppose that $n_j$ is prime to $p$. Then the
    specialization $\xb_j$ of $x_j$ lies on a primitive tail $\Xb_v$,
    and the ramification index of $\Yb_v\to\Xb_v$ at $\xb_j$ is $n_j$.
  \item
    With the same notation as in (3), suppose that $n_j=pn_j'$. Then
    $\xb_j$ lies on a component $\Xb_v$ which is not a tail. Moreover,
    the Galois cover $\Zb_v\to\Xb_v$ has ramification index
    $n_j'$ at $\xb_j$. 
  \end{enumerate}
\end{lem}

\begin{proof} This is essentially proved in \cite{Raynaud98}, \S 2 and \S 3.
There are two points where our assumptions differ from the assumptions
made in \cite{Raynaud98}. Firstly, it is assumed in \cite{Raynaud98}
that $p$ exactly divides the order of $G$. However, one only needs
that $f$ has mild reduction. The second point is that (3) becomes
false if we replace our definition of the stable reduction with the
definition used in \cite{Raynaud98}. It is proved in \cite{modular},
Lemma 2.1.1 that (3) holds if we define $\fb$ as in Definition
\ref{stabledef}.  \end{proof}

\subsection{Deformation data} \label{defdat}

\subsubsection{}  \label{defdat1}

A priori, our definition of a deformation datum is independent of the
study of stable reduction.  So the following notation replaces, during
\S \ref{defdat1}, the notation introduced so far. Let $p$ be a prime
number, $k$ an algebraically closed field of characteristic $p$ and
$X$ a smooth projective and irreducible curve over $k$.

\begin{defn} \label{defdatdef}
  A {\em deformation datum} over $X$ is a pair $(Z,\omega)$, where
  $Z\to X$ is a tame (connected) Galois cover and $\omega$ is a
  meromorphic differential form on $Z$ such that the following holds.
  \begin{enumerate}
  \item
    The order of the Galois group $H=\Gal(Z/X)$ is prime to $p$. 
  \item
    We have
    \[
        \beta^*\omega \;=\; \chi(\beta)\cdot\omega, \qquad
           \text{for all $\beta\in H$,}
    \]
    where $\chi:H\to\FF_p^\times$ is a character with values in the
    prime field $\FF_p$. 
  \item 
    The differential $\omega$ is either logarithmic (i.e.\ of the
    form $\omega=\diff u/u$) or exact (i.e.\ $\omega=\diff u$). In the
    first case, the deformation datum $(Z,\omega)$ is called {\em
      multiplicative}, in the second case it is called {\em additive}.
  \end{enumerate}
  The pair $(H,\chi)$ is called the {\em type} of the deformation
  datum $(Z,\omega)$. 
\end{defn}

Let $\xi\in Z$ be a closed point and $\tau\in X$ its image. Denote the
stabilizer of $\xi$ by $H_\xi\subset H$. We attach to the triple
$(Z,\omega,\tau)$ the following three invariants:
\begin{equation} 
    m_\tau \;:=\; |H_\xi|, \qquad h_\tau \;:=\; \ord_\xi(\omega)+1, \qquad
       \sigma_\tau \;:=\; h_\tau/m_\tau.
\end{equation}
It is clear that $m_\tau$, $h_\tau$ and $\sigma_\tau$ depend only on
$\tau$ but not on $\xi$, and that $(m_\tau,h_\tau)=(1,1)$ for all but
a finite number of points $\tau\in X$.  If $(m_\tau,h_\tau)\not=(1,1)$
then $\tau$ is called a {\em critical point} for the deformation datum
$(Z,\omega)$.

Given a deformation datum $(Z,\omega)$, we often use the following
or a similar notation. Let $(\tau_j)_{j\in J}$ be the family of critical
points for $(Z,\omega)$, indexed by the finite set $J$.  For $j\in J$,
we write $m_j$, $h_j$ and $\sigma_j$ instead of $m_{\tau_j}$, $h_{\tau_j}$
and $\sigma_{\tau_j}$. The data $(m_j;h_j)_{j\in B}$ is called the {\em
  signature} of the deformation datum $(Z,\omega)$. A straightforward
computation shows that
\begin{equation} \label{defdateq2}
   \sum_{j\in J}\, (\sigma_j-1) \;=\; 2g_X-2.
\end{equation}
(where $g_X$ denotes the genus of the curve $X$). 

\begin{rem} \label{defdatrem}
  Let $\tilde{Z}:=Z/\Ker(\chi)$. It is
  easy to see that $\omega$ descents to a differential
  $\tilde{\omega}$ on $\tilde{Z}$ and that
  $(\tilde{Z},\tilde{\omega})$ is again a deformation datum, with
  critical points $\tau_j$ and invariants 
  \[
      \tilde{m}_j \;=\; \frac{m_j}{|\Ker(\chi)|}, \qquad
      \tilde{h}_j \;=\; \frac{h_j}{|\Ker(\chi)|}, \qquad
      \tilde{\sigma}_j \;=\; \frac{\tilde{h}_j}{\tilde{m}_j} 
           \;=\; \sigma_j.
  \]
\end{rem}

\subsubsection{} \label{defdat2}

Let us now return to our original notation, concerning the Galois
cover $f:Y\to X$ and its stable reduction $\fb:\Yb\to\Xb$. Let $\Xb_v$
be an interior component of $\Xb$ (see Definition
\ref{taildef}). Following \cite{YannickArbres} and
\cite{SaidiTorsors}, we define a deformation datum $(\Zb_v,\omega_v)$
over $\Xb_v$, as follows.

Choose an irreducible component $\Yb_v$ of $\Yb$ above $\Xb_v$. By
Lemma \ref{stablem} (2), the inertia group $I_v$ of $\Yb_v$ is cyclic
of order $p$ and the quotient group $H_v:=G_v/I_v$ is of prime-to-$p$
order. Moreover, the cover $\Yb_v\to\Xb_v$ is the composition of a
purely inseparable cover $\Yb_v\to\Zb_v$ of degree $p$ and an
$H_v$-Galois cover $\Zb_v\to\Xb_v$. Let $\Yd_v:=\Spec\Od_{Y_R,v}$ be
the completion of $Y_R$ at the generic point of the component $\Yb_v$;
thus $\Od_{Y_R,v}$ is a complete discrete valuation ring with residue
field $k(\Yb_v)$. The action of $I_v$ on $Y_R$ induces an action on
$\Yd_v$.  Let $\Zd_v:=\Yd_v/I_v$ be the quotient; we have $\Zd_v=\Spec
A$ for a complete discrete valuation ring $A$ with residue field
$k(\Zb_v)$. We may assume that the field $K$ contains the $p$th roots
of unity.  Let us choose an isomorphism $I_v\cong\bmu_p(K)$ between $I_v$
and the group of $p$th roots of unity. By \cite{Raynaud98}, \S 1.2.1,
the map $\Yd_v\to\Zd_v$ is a torsor under a finite flat $R$-group
scheme $\G$ of rank $p$, such that $\G\otimes_R K=\bmu_p$ acts on
$\Yd_v\otimes_R K$ via the isomorphism $I_v\cong\bmu_p(K)$. By Kummer
theory, the $\bmu_p$-torsor $\Yd_v\otimes_R K\to\Zd_v\otimes_R K$ is
given by a Kummer equation $y^p=u$, with $u\in A$.

To define the differential $\omega_v$, we have to distinguish two
cases. First, let us assume that $\G\otimes_R k\cong\bmu_p$. This
happens if and only if we may choose the Kummer equation $y^p=u$ such
that the residue $\ub$ of $u$ is not a $p$th power. In this case, we
set $\omega_v:=\diff\ub/\ub$ and we say that the cover $f$ has {\em
multiplicative reduction} at $\Yb_v$. Otherwise $\G\otimes_R
k\cong\balpha_p$ and we say that $f$ has {\em additive reduction} at
$\Yb_v$. In this case we can choose the equation $y^p=u$ such that
$u=1+\pi^pw$, with $\pi\in R$, $w\in A$ and $0<v_R(\pi)<v_R(p)/(p-1)$,
and such that the residue $\wb$ of $w$ is not a $p$th power. We set
$\omega_v:=\diff\wb$. It is clear that $\omega_v$ depends on the
choice of the isomorphism $I_v\cong\bmu_p(K)$. Indeed, composing this
isomorphism with the automorphism $\zeta_p\mapsto\zeta_p^a$ (where
$a\in\FF_p^\times$) amounts to multiplying $\omega_v$ by $a$. Let
$\chi_v:H_v\to\FF_p^\times$ be the character describing the action of
$H_v$ on $I_v$ by conjugation. It is clear that
$\beta^*\omega_v=\chi_v(\beta)\cdot\omega_v$ for all $\beta\in H_v$.

We have defined the differential $\omega_v$ over the generic point of
$\Zb_v$. However, we can extend this definition to a neighborhood of
any point of $\Zb$ which lies over a smooth point of $\Xb$. To be more
precise, let $\tau\in\Xb_v$ be a closed point which is smooth in $\Xb$
and choose points $\eta\in\Yb_v$ and $\xi\in\Zb_v$ above $\tau$. Let
$\Yd_\eta$ be the completion of $Y_R$ at $\eta$ and
$\Zd_\xi=\Yd_\eta/I_v$. Then
$\Yd_\eta\otimes_Rk=\Spec\Od_{\Yb_v,\eta}\cong k[[y]]$ and
$\Zd_\xi=\Spec\Od_{\Zb_v,\xi}\cong k[[z]]$. As above, the action of
$I_v$ on $\Yd_\eta$ gives rise to a differential $\omega_\xi$ of the
field ${\rm Frac}(\Od_{\Zb_v,\xi})\cong k((z))$. One checks easily
that $\omega_\xi$ is the pullback of $\omega_v$ to the completion of
$\Zb_v$ at $\xi$. Now the local definition of $\omega_\xi$ shows the
following. The differential $\omega_v$ does not vanish at $\xi$.
Moreover, if $\omega_v$ has a pole at $\xi$, then this pole is simple,
$\omega_\xi$ is logarithmic and the point $\tau\in\Xb_v$ is the
specialization of a branch point of $f$ with ramification index
divisible by $p$. We summarize the discussion in the next proposition.

\begin{prop} \label{defdatprop1}
  Over each interior component $\Xb_v$ we obtain a deformation datum
  $(\Zb_v,\omega_v)$, with the following properties.  If
  $\tau\in\Xb_v$ is a critical point for $(\Zb_v,\omega_v)$ then one
  of the following two statements holds.
  \begin{enumerate}
  \item
    The point $\tau\in\Xb_v$ is a singular point of $\Xb$.
  \item
    The point $\tau$ is the
    specialization of a point $x_j\in S_0$ at which the cover $f$ has
    ramification of order divisible by $p$. 
  \end{enumerate}
  In Case (1), we have $\sigma_{\tau}\not\in\{0,1\}$. In Case (2), we
  have $\sigma_{\tau}=0$ and $(\Zb_v,\omega_v)$ is multiplicative.
\end{prop}

\subsubsection{} \label{defdat3}

By Proposition \ref{defdatprop1}, the stable reduction $\fb:\Yb\to\Xb$
of a $G$-cover $f:Y\to X$ with mild reduction is naturally equipped
with a collection of deformation data $(\Zb_v,\omega_v)$. These
deformation data are related among each other by certain compatibility
conditions. To formulate these conditions, let $\tau$ be a singular
point of $\Xb$, and let $\Xb_{v_1},\Xb_{v_2}$ be the two components of
$\Xb$ which intersect in $\tau$. For $i=1,2$, we define a pair
$(m_i,h_i)$ as follows.
  \begin{itemize}
  \item If $\Xb_{v_i}$ is an interior component, we define $(m_i,h_i)$
    as the pair of invariants attached to
    $(\Zb_{v_i},\omega_{v_i},\tau)$.
  \item Suppose that $\Xb_{v_i}=\Xb_j$ is a tail of $\Xb$. Let
    $\eta\in\Yb_j$ be a point lying above $\tau$ and let $I_{\eta}\subset
    G$ be the inertia group at $\eta$. We define $m_i$ as the order of
    the maximal prime-to-$p$ quotient of $I_{\eta}$. Furthermore, we
    let $h_i$ be the conductor of $I_{\eta}$, i.e.\ the unique integer
    such that $|I_{\eta,h_i}|=p$ and $I_{\eta,h_i+1}=\{1\}$. (Here
    $(I_{\eta,\nu})_{\nu\geq 0}$ denotes the filtration of higher
    ramification groups.)
  \end{itemize}

\begin{prop} \label{defdatprop2}  
  We have $m_1=m_2$ and $h_1=-h_2$.
\end{prop} 

\begin{proof} The equality $m_1=m_2$ is trivial. Let
$\Yd:=\Spec\Od_{Y_R,\eta}$ be the completion of $Y_R$ at $\eta$. The
$R$-scheme $\Yd$ is a formal annulus and is acted upon by $I_{\eta}$.
The invariants $h_1$ and $h_2$ can be defined purely in terms of this
local action, and the equality $h_1=-h_2$ is proved in
\cite{YannickArbres}, Proposition 1.10. (In \cite{YannickArbres}, the
invariant $h$ is called $-m$.)  \end{proof}

\begin{rem}
  It seems plausible that the existence of a set of deformation data
  satisfying these compatibility conditions gives a necessary and
  sufficient condition for a map $\fb:\Yb\to\Xb$ between semistable
  curves to occur as the stable reduction of a Galois cover between
  smooth curves in characteristic $0$. That this is true for three
  point covers is proved in \S \ref{main}. In the case of $p$-cyclic
  covers, results in this direction are obtained in
  \cite{YannickArbres} and \cite{SaidiTorsors}.
\end{rem}

\subsection{The vanishing cycle formula} \label{vcf}

The vanishing cycle formula, proved by Raynaud in \cite{Raynaud98}, is
the key formula which controls the (local) invariants describing the
bad reduction occuring over the tails by a global expression.  In this
subsection we give another proof. The main observation is that
\eqref{defdateq2} gives a `local' vanishing cycle formula for each
interior component of $\Xb$. The global vanishing cycle formula is
then easily obtained from its local version using the compatibility
relation of Proposition \ref{defdatprop2} and induction over the tree
of components of $\Xb$. We use this occasion to set up some more
notation which may seem at bit technical at first but will pay off in
\S \ref{simple}.

\subsubsection{} \label{vcf1}

Recall that an {\em unordered graph} is given by a quintupel
$(V,E,s,t,\bar{\;\;})$. Here $V$ is the set of {\em vertices}, $E$ is
the set of {\em (ordered) edges}, $s,t:E\to V$ are the {\em source}
and the {\em target} map, and $E\iso E$, $e\mapsto\bar{e}$, is a fixed
point free involution of $E$ such that $s(\bar{e})=t(e)$ for all $e\in
E$.

Let $T'=(V',E',s,t,\bar{\;\;})$ be the dual graph of the semistable
curve $\Xb$. Thus, a vertex $v\in V'$ corresponds to an irreducible
component $\Xb_v$ of $\Xb$, and an edge $e\in E'$ corresponds to a
triple $(\xb,v_1,v_2)$ such that the (distinct) components $\Xb_{v_1}$
and $\Xb_{v_2}$ intersect in the (singular) point $\xb$. Let
$B\wild\subset B_0$ be the set of indices $j$ such that ramification
index of the Galois cover $f:Y\to X$ at $x_j$ is divisible by $p$.  We
define $T=(V,E,s,t,\bar{\;\;})$ as the unordered graph obtained from
$T'$ by adding a vertex for each point $x_j$ with $j\in B\wild$. More
precisely, $V:=V'\cup B\wild$ and a vertex $j\in B\wild$ is connected
to the vertex $v\in V'$ by an edge $e\in E$ if and only if the point
$x_j\in S_0$ specializes to the component $\Xb_v$. Note that $T$ is a
connected tree. Also, to each edge $e$ of $T$ corresponds a closed
point $\xb_e\in\Xb$, as follows.  If the edge $e$ joins two vertices
$v_2,v_2\in V'$ then $\xb_e$ is the point where $\Xb_{v_1}$ intersects
$\Xb_{v_2}$. If $e$ joins two vertices $j\in B\wild$ and $v\in V'$
then $\xb_e:=\xb_j\in\Xb_v$ is the point to which $x_j$ specializes.

We consider the vertex $v_0$ corresponding to the original component
$\Xb_{v_0}$ as the {\em root} of the tree $T$. Also, we denote by $B$
the set of {\em leaves} of $T$, i.e.\ the set of vertices $v\in
V-\{v_0\}$ which are adjacent to a unique vertex $v'$. Thus, a vertex
$v\in V\intern:=V-B$ corresponds to an interior component $\Xb_v$. By
definition, $B\wild\subset B$. An element $j$ in $B\prim=B_0-B\wild$
corresponds to a ramification point $x_j$ of $f$ of order
prime-to-$p$. By Lemma \ref{stablem}, $x_j$ specializes to a point
$\xb_j$ on a primitive tail $\Xb_v$. The resulting map $B\prim\to B$
is clearly injective. Therefore, we may and will consider
$B_0=B\prim\cup B\wild$ as a subset of $B$. Then $B\new:=B-B_0$ is the
set of leaves corresponding to the new tails of $\Xb$.

\subsubsection{} \label{vcf2}

We shall now attach integers $m_e,h_e$ and a rational number
$\sigma_e$ to each edge $e$ of the tree $T$. Suppose first that
$v:=s(e)$ is an interior vertex. Let $(\Zb_v,\omega_v)$ be the
deformation datum over $\Xb_v$, as defined in \S \ref{defdat2}. Note
that $\xb_e\in\Xb_v$ is a critical point for $(\Zb_v,\omega_v)$, by
Proposition \ref{defdatprop1}. We define the numbers $m_e,h_e$ and
$\sigma_e$ as the invariants attached to the triple
$(\Zb_v,\omega_v,\xb_e)$. If $v:=s(e)$ is a leaf of $T$, then we set
$m_e:=m_{\bar{e}}$, $h_e:=-h_{\bar{e}}$ and
$\sigma_e:=-\sigma_{\bar{e}}$. Equation \eqref{defdateq2} and
Proposition \ref{defdatprop2} imply the following.

\begin{prop} \label{vcfprop}
  For all interior vertices $v$ of $T$ we have
  \begin{equation} \label{vcfpropeq1}
       \sum_{s(e)=v}\, (\sigma_e-1) \;=\; 2g_v - 2
  \end{equation}
  (where $g_v$ is the genus of $\Xb_v$). Furthermore, for all edges
  $e$ of $T$ we have
  \begin{equation} \label{vcfpropeq2}
        \sigma_e + \sigma_{\bar{e}} \;=\; 0.
  \end{equation}
\end{prop}

For each leaf $j\in B$, let $e_j$ be the unique edge with target $j$.
We set $m_j:=m_{e_j}$, $h_j:=m_{e_j}$ and $\sigma_j:=\sigma_{e_j}$. If
$j\not\in B\wild$ then $m_j$ is the prime-to-$p$ part of the
ramification index of the Galois cover $\Yb_j\to\Xb_j$ at
$\xb_e=\infty_j$ and $h_j$ is the conductor of the wild ramification
in that point, see Proposition \ref{defdatprop2}. Therefore, the
rational number $\sigma_j=h_j/m_j$ is the ramification invariant
introduced in \cite{Raynaud98}, \S 1.1. It is equal to the jump in the
filtration of higher ramification groups for the branch point
$\infty_j\in\Xb_j$, with respect to the upper numbering, see
\cite{SerreCL}. Note that we have defined $\sigma_j=0$ for $j\in
B\wild$. We shall call the tuple $(\sigma_j)_{j\in B}$ the {\em
signature} of the stable reduction $\fb:\Yb\to\Xb$. It is shown in
\cite{Raynaud98} that $\sigma_j>1$ for $j\in B\new$. Moreover, we have
the following.

\begin{cor}[The vanishing cycle formula] \label{vcfcor}
  The signature $(\sigma_j)_{j\in B}$ of $\fb$ verifies
  \begin{equation} \label{vcfeq}
     \sum_{j\in B\prim}\,\sigma_j \;+\; 
        \sum_{j\in B\new}\,(\sigma_j-1) \;=\; 2g_X -2 + |B_0|.
  \end{equation}
\end{cor}

\begin{proof} This is proved in \cite{Raynaud98}, \S 3.4. Here is another
proof, using Proposition \ref{vcfprop}. Let $T_0\subset
T_1\subset\ldots\subset T_n$ be a chain of subtrees of $T$ such that
$T_0=\{v_0\}$, $T_i=T_{i-1}\cup\{v_i\}$ and $T_n=T-B$. Let $E_i$
denote the set of edges of $T$ such that $s(e)\in T_i$ and
$t(e)\not\in T_i$. We claim that $\sum_{e\in E_i} (\sigma_e-1)=2g_X-2$
holds for all $0\leq i\leq n$.  Indeed, Proposition \ref{vcfprop}
shows that the claim is true for $i=0$ and then, by induction, for all
$i$. For $i=n$, the claim is equivalent to \eqref{vcfeq}.  \end{proof}


\section{The case of three point covers} \label{threepoints}

In this section we apply the results of the last section to the case
of three point covers. In \S \ref{simple} we prove our first main
result, which essentially states that the stable reduction of a three
point cover with mild reduction is as simple as possible. The basic
idea underlying the proof of this result is that $\fb:\Yb\to\Xb$,
together with the associated deformation data, should be a `rigid'
objects, i.e.\ should not allow any nontrivial deformation. This
imposes strong restrictions on the deformation data $(\Zb_v,\omega_v)$
attached to $\fb$. We verify these restrictions coming from rigidity
only indirectly, using the local and global vanishing cycle
formulas. A key lemma from \cite{special} shows that $\fb$ can support
no additive deformation datum, and only a single multiplicative
one. To make this proof work one has to exclude two exceptional cases,
see \S \ref{except}. This is somewhat annoying but not serious,
because all our final results apply to these exceptional cases as
well.

In \S \ref{Gdefdat} we define the notion of a {\em special
  $G$-deformation datum}. To such an object we can associate a map
$\fb:\Yb\to\Xb$ which looks as if it were the stable reduction of a
three point $G$-cover, and which we call a {\em special $G$-map}.

\subsection{The structure of $\fb$} \label{simple}

\subsubsection{} \label{except}

We continue with the notation introduced in \S \ref{stablered}. In
addition, we assume that $X_0=\PP^1_{R_0}$ and
$S_0=\{0,1,\infty\}$. Note that the right hand side of \eqref{vcfeq}
is equal to $1$.

\begin{defn} \label{exceptdef}
  We say that the stable reduction of $f$ is {\em exceptional} if one
  of the following two cases occurs.
  \begin{enumerate}
  \item
     We have $B\prim=\{j_0\}$, $B\new=\emptyset$ and $\sigma_{j_0}=1$.
  \item
     We have $B\prim=\emptyset$, $B\new=\{j_0\}$ and $\sigma_{j_0}=2$.
  \end{enumerate}
\end{defn}

\begin{prop} \label{exceptprop}
  Let $f:Y\to X=\PP^1_{\Kb}$ be a Galois cover ramified at
  $\{0,1,\infty\}$, with exceptional reduction.  Then the
  curve $\Xb$ consists only of the original component $\Xb_0$ and the
  unique tail $\Xb_{j_0}$.
\end{prop}

\begin{proof} We will prove the proposition only for Case (2) of Definition
\ref{exceptdef} and leave Case (1) to the reader. Let $T'$ be the dual
graph of the curve $\Xb$, considered as a subtree of $T$. By
assumption, $T'$ has a unique leaf $j_0$. Therefore, $T'$ is simply a
chain going from the root $v_0$ to $j_0$. Let $v$ be the unique vertex
of $T'$ adjacent to the leaf $j_0$ and let $e$ be the edge with source
$v$ and target $j_0$. This means that the interior component $\Xb_v$
intersects the unique (new) tail $\Xb_{j_0}$ in the point $\xb_e$. We
have to show that $v=v_0$, i.e.\ that $\Xb_v$ is the original
component $\Xb_0$. So, let us assume that $v\not=v_0$. Then there
exists a unique edge $e'\not=e$ with source $v$, corresponding to a
point $\xb_{e'}$ where $\Xb_v$ intersects another component
$\Xb_{v'}$. Let us first suppose that there exists another edge $e''$
(of the tree $T$) with source $v$ besides $e$ and $e'$. Then
$t(e'')\in B\wild$, i.e. the point $\xb_{e''}$ is the specialization
of a branch point of $f$ with ramification index divisible by $p$.
Since $\Xb_v$ is not the original component, there can be at most one
branch point specializing to $\Xb_v$. Therefore, $e$, $e'$, $e''$ are
the only edges with source $v$. Since $\sigma_e=2$ and
$\sigma_{e''}=0$, \eqref{vcfpropeq1} shows that $\sigma_{e'}=-1$. This
means that the differential $\omega_v$ has a pole of order $\geq 2$ at
each point of $\Zb_v$ lying above $\xb_{e'}$. This shows that
$\omega_v$ cannot be logarithmic. On the other hand, $\omega_v$ has a
simple pole at each point above $\xb_{e''}$, so $\omega_v$ cannot be
exact. This is a contradiction. We may therefore assume that $e$ and
$e'$ are the only edges with source $v$. It follows that the cover
$\Zb_v\to\Xb_v$ is ramified at most at two points; in particular, the
genus of $\Zb_v$ and therefore the genus of $\Yb_v$ as well, is zero.
Moreover, $\Yb_v$ intersects the rest of $\Yb$ in two points and no
ramification point specializes to $\Yb_v$. Thus, the `three point
condition' does not hold for $\Yb$, in contradiction with the
definition of the stable model in \S \ref{stable}. The proposition
follows.  \end{proof}

\begin{rem} \label{goodredrem}
  It may happen that the Galois cover $f:Y\to X$ has bad reduction but
  that the curve $Y$ has good reduction. If this is the case then the
  stable reduction of $f$ is exceptional.
\end{rem}

\subsubsection{}  \label{combi}

For the rest of \S \ref{simple} we assume that the reduction of
$f$ is not exceptional. The vanishing cycle formula \eqref{vcfeq},
together with the inequalities $\sigma_j>0$ for $j\in B\prim$ and
$\sigma_j>1$ for $j\in B\new$, immediately imply the inequalities
\begin{equation} \label{sigmaineq}
  0<\sigma_j<1 \quad\text{for $j\in B\prim$,}\qquad
  1<\sigma_j<2 \quad\text{for $j\in B\new$.}
\end{equation}
Moreover, $|B\prim\cup B\new|\geq 2$.  For an edge $e$ of $T$ we
define $\nu_e:=\lfloor\,\sigma_e\,\rfloor$ to be the largest integer
less than or equal to the rational number $\sigma_e$. Note that we
have $\sigma_e=\nu_e+\gen{\sigma_e}$, where $\gen{\sigma_e}\geq 0$
denotes the fractional part of $\sigma_e$. For a leaf $j\in B$, we
write $\nu_j:=\lfloor\,\sigma_j\,\rfloor$. By \eqref{sigmaineq} we
have $\nu_j=0$ for $j\in B_0$ and $\nu_j=1$ for $j\in B\new$. Since
$|B_0|=3$, the vanishing cycle formula \eqref{vcfeq} implies
\begin{equation} \label{combieq3}
   \sum_{j\in B}\; \gen{\sigma_j} \;=\; 1.
\end{equation}
Recall that $T'$ denotes the dual graph of the curve $\Xb$, considered
as subtree of $T$. Let $E'$ denote the set of edges of $T'$. Moreover,
let $E'\intern\subset E'$ be the set of edges connecting two interior
vertices.  

\begin{prop} \label{aeprop}\ 
\begin{enumerate}
\item For all edges $e\in E'$ we have $\nu_e+\nu_{\bar{e}}=-1$.
\item For all interior vertices $v\in V\intern$ we have
  \[
        \sum_{s(e)=v}\; \gen{\sigma_e} \;=\; 1 \qquad\text{and}\qquad
         \sum_{s(e)=v}\,(\nu_e-1) \;=\; -3.
  \]
\end{enumerate}
\end{prop}

\begin{proof} 
For $e\in E'$ define 
\[
     d_e \;:=\; \gen{\sigma_e}+\gen{\sigma_{\bar{e}}}-1 \;=\;
       -\nu_e-\nu_{\bar{e}}-1.
\]
The last equality follows from \eqref{vcfpropeq2}. Clearly,
$d_e=d_{\bar{e}}$, $d_e\in\{-1,0\}$ and $d_e=-1$ if and only if
$\sigma_e\in\ZZ$.  For $v\in V\intern$ define
\[
     c_v \;:=\; -1 + \sum_{s(e)=v}\gen{\sigma_e} \;=\;
        -3 + \sum_{s(e)=v} (\nu_e-1).
\]
The last equality is a consequence of \eqref{vcfpropeq1}. It is easy
to see that $c_v\in\{-1,0,1,\ldots\}$ and that $c_v=-1$ if and only if
$d_e=-1$ for all $e\in E'$ such that $s(e)=v$. We have to show that
$d_e=c_v=0$.

It is convenient to introduce the following terminology. A vertex $v$
of $T$ is said to {\em precede} an edge $e$ (notation: $v\prec e$) if $v$
lies in the connected component of $T-\{e\}$ which contains $s(e)$.
An edge $e'$ is said to precede another edge $e$ (notation: $e'\prec
e$) if $t(e')\prec e$. Note that `$\prec$' is a partial order on the
set $E$ of edges. 

Suppose that there exists an edge $e\in E'$ with $d_e\not=0$. We may
suppose that $e$ is minimal with this property, i.e.\ $d_{e'}=0$ for
all edges $e'\in E'$ which precede $e$. We claim that $c_v\geq 0$ for
all $v\in V\intern$ which preceed $e$. Indeed, since $\Xb$ has at least
two tails, there exist at least two (distinct) edges in $E'$
with target $v$. At least one of them precedes $e$. Hence there exists
an edge $e'\in E'$ with $d_{e'}=0$ and $t(e)=v$. This proves the
claim. As a consequence, we obtain the inequality
\begin{equation} \label{aepropeq1}
\begin{split}
    0 \;\leq\; \sum_{\begin{array}{c} v\in V\intern \\ 
                      v\prec e \end{array}} c_v 
   &\;=\;\gen{\sigma_e} \;+\; 
      (\sum_{\begin{array}{c} e'\in E'\intern \\ 
                     e'\prec e \end{array}}\, d_{e'}\,) 
  \,+\,(\sum_{\begin{array}{c} j\in B \\ j\prec e \end{array}} 
          \gen{\sigma_j}\,) \,-\,1 \\
   & \;=\; \gen{\sigma_e} \;+\; 
       (\sum_{\begin{array}{c}  j\in B \\ j\prec e \end{array}} 
            \gen{\sigma_j}\,) \,-\, 1.
\end{split}
\end{equation}
Since the curve $\Xb$ has at least two tails, there exists $j\in
B\prim\cup B\new$ such that $j$ does not precede $e$. Moreover
$\sigma_j>0$. Therefore, \eqref{aepropeq1} implies $\gen{\sigma_e}>0$.
But this is a contradiction to the assumption $d_e\not=0$. We conclude
that $d_e=0$ for all $e\in E'$.  It remains to prove $c_v=0$ for all
$v\in V\intern$. The same argument we used above together with $d_e=0$
for all $e\in E'$ shows that $c_v\geq 0$ for all $v\in V\intern$. Now
the proposition follows from
\[
    \sum_{v\in V\intern}\, c_v \;=\; 
 (\,\sum_{e\in E\intern'}\, d_e\,) \;+\; (\,\sum_{j\in B}\,\gen{\sigma_j}\,) 
        -1 \;=\; 0.
\]
\end{proof}

\begin{lem} \label{nuelem}
  Suppose that $\fb$ is not exceptional, and let $e$ be an edge of $T$.
  Then 
  \[
      -2 \leq \nu_e \leq 1.
  \]
  Moreover, $\nu_e\geq 0$ if and only if the root $v_0$ precedes
  $e$. In particular, if $v$ is an interior vertex unequal to $v_0$
  then there exists a unique edge $e$ with source $v$ such that
  $\nu_e<0$.
\end{lem}

\begin{proof} (Compare to \cite{special}, Lemma 2.6.)  Using Proposition
\ref{aeprop} and induction in the tree $T$ one shows that
\[
     \nu_e \;=\; 1 - |\,\{\,j\in B_0\mid j>e \,\}\,|.
\]
Note that $|B_0|=3$ and that the root $v_0$ is the median of the set
of leaves $B_0$. The lemma follows immediately. 
\end{proof}

\subsubsection{}

We can now prove our first main result. It applies in particular if
$p$ strictly divides the order of $G$. 

\begin{thm} \label{thm1}
  Let $f:Y\to X=\PP^1_{\Kb}$ be a Galois cover ramified at
  $\{0,1,\infty\}$, with bad reduction. If $f$ has mild reduction
  (Definition \ref{milddef}) then the curve $\Xb$ consists only of the
  original component $\Xb_0$ and the tails $\Xb_j$, $j\in B\prim\cup
  B\new$. 
\end{thm}

\begin{proof} If $B\wild=\emptyset$ then this follows from \cite{special},
Theorem 2.1, using Raynaud's construction of the {\em auxiliary
  cover}.  The proof we give below follows the same lines, but avoids
the construction of the auxiliary cover. 

We prove the theorem by contradiction. Assume that there exists an
interior vertex $v$ of $T$ other than the root $v_0$. By Proposition
\ref{exceptprop}, the reduction of $f$ is not exceptional.  We may
assume that $v$ is adjacent to the root $v_0$, i.e.\ there exists an
edge $e_0$ with source $v$ and target $v_0$. This means that the
component $\Xb_v$ intersects the original component $\Xb_0$ in the
point $\xb_{e_0}$.  Let $(\Zb_v,\omega_v)$ be the deformation datum
over $\Xb_v$ as defined in \S \ref{defdat2}. For the rest of this
proof, $e$ runs over the set of edges with source $v$.  By Proposition
\ref{aeprop} we have $\sum_e \gen{\sigma_e}=1$ and $\sum_e
(\nu_e-1)=-3$. Also, Lemma \ref{nuelem} implies that $-2\leq\nu_e\leq
1$ for all $e$ and $\nu_e<0$ if and only if $e=e_0$. In this
situation, \cite{special}, Lemma 2.8 says that $\omega_v$ cannot be
exact (actually, in \cite{special} one assumes that the type of
$(\Zb_v,\omega_v)$ has an injective character $\chi_v$; the general
case follows from Remark \ref{defdatrem}).  Therefore, $\omega_v$ is
logarithmic. On the other hand, since $\nu_{e_0}<0$, the differential
has a pole of order $\geq 2$ at each point of $\Zb_v$ lying above
$\xb_{e_0}$. Therefore, $\omega_v$ cannot be logarithmic. This gives
the desired contradiction and proves the theorem.  \end{proof}

\subsection{Special $G$-deformation data} \label{Gdefdat}

In this section we define special $G$-deformation data. To a special
$G$-deformation datum we then associate a special $G$-map
$\fb:\Yb\to\Xb$. Briefly, a special $G$-map is a map between
semistable curves over $k$ which looks like the stable reduction of a
three point $G$-cover with mild reduction.

\subsubsection{} \label{Gdefdat1}

Let $k$ be an algebraically closed field of characteristic $p>0$ and 
$\Xb$ a smooth projective curve over $k$ of genus $0$. Let
$(\Zb,\omega)$ be a deformation datum over $\Xb$, of type
$(H,\chi)$. As usual, we denote by $(\tau_j)_{j\in B}$ the set of
critical points for $(\Zb,\omega)$. Also, for $j\in B$ we denote by
$m_j$, $h_j$ and $\sigma_j$ the invariants attached to
$(\Zb,\omega,\tau_j)$. See \S \ref{defdat}.

\begin{defn} \label{specialdefdatdef}
  The deformation datum $(\Zb,\omega)$ is called {\em special} if for
  all $j\in B$ we have 
  $\sigma_j<2$ and $\sigma_j\not=1$ and if there exists a subset
  $B_0\subset B$ with exacly three elements such that $\sigma_j<1$ if
  and only if $j\in B_0$. 
\end{defn}

\begin{rem} \label{specialrem}
  In case the character $\chi$ is injective and $\sigma_j>0$ we
  recover the definition of a special deformation datum given in
  \cite{special}.
\end{rem}

In the following we assume that $(\Zb,\omega)$ is special. Let
$B\wild$ (resp.\ $B\prim$, resp.\ $B\new$) denote the set of all
indices $j\in B$ such that $\sigma_j=0$ (resp.\ $0<\sigma_j<1$, resp.\ 
$1<\sigma_j<2$). Then $B_0=B\wild\cup B\prim$, and it follows
from \eqref{vcfpropeq1} that 
\begin{equation} \label{Gdefdateq1}
      \sum_{j\in B}\; \gen{\sigma_j} \;=\; 1.
\end{equation}
We call the special deformation datum $(\Zb,\omega)$ {\em normalized}
if $\Xb=\PP^1_k$ and $\{\tau_j\mid j\in B_0\}=\{0,1,\infty\}$.

Their connection to three point covers suggests that special
deformation data are `rigid' objects, i.e.\ there should not
exist any nontrivial families of special deformation data. This is
indeed true; a precise statement is formulated and proved in
\cite{cotang}, \S 5.4. See also \S \ref{aux} below. 





\subsubsection{}  \label{Gdefdat2}

In this subsection, $\Xb$ will denote a smooth projective curve over
$k$ of genus $0$, together with a distinguished closed point
$\infb\in\Xb$. We also fix a finite group $G$.

\begin{defn} \label{tailcoverdef}
  A {\em $G$-tail cover} is a (not necessarily connected) $G$-cover
  $\fb:\Yb\to\Xb$ such that the following holds.
  \begin{enumerate}
  \item
    The cover $\fb$ is wildly ramified above $\infb$, of order $pm$ with
    $(p,m)=1$.
  \item
    The restriction of $\fb$ to $\Xb-\{\infb\}$ is tame and branched 
    at most at one point.
  \end{enumerate}
  A {\em pointed $G$-tail cover} is a $G$-tail cover together with a
  point $\eta\in\Yb$ above $\infb$ and a subgroup $H$ of the inertia
  group $I_\eta$ such that $I_\eta=I_{\eta,1}\rtimes H$. 
  An (outer) automorphism of a pointed $G$-tail cover is an
  automorphism $\gamma:\Yb\iso\Yb$ which normalizes the action of $G$
  (and hence induces an automorphism of $\Xb$), centralizes the action
  of $H\subset G$ and fixes the point $\eta$. Such an automorphism is
  called {\em inner} if it centralizes the action of $G$.
\end{defn}

Let $(\fb:\Yb\to\Xb,\eta,H)$ be a pointed $G$-tail cover. By Condition
(1), $I_{\eta,1}$ is cyclic of order $p$ and $H$ cyclic of order
prime-to-$p$. Let $m$ denote the order of $H$ and $h$ the conductor of
$I_{\eta,1}$. The pair $(m_j,h_j)$ is called the {\em ramification
type} and the rational number $\sigma:=h/m$ the {\em ramification
invariant} of $\fb$. Note that $h>0$ and $(h,p)=1$.  Moreover,
$m|(p-1)h$, by the Hasse--Arf Theorem (see \cite{SerreCL}, IV \S 3).

\begin{defn}  \label{specialtaildef}
  A tail cover $\fb:\Yb\to\Xb$ is called {\em special} if
  the following two conditions hold.
  \begin{enumerate}
  \item
    We have $\sigma<2$ and $\sigma\not=1$.
  \item
    If $\sigma>1$ then the restriction of $\fb$ to $\Xb-\{\infb\}$
    is \'etale. 
  \end{enumerate}
  If $\sigma<1$ then we say that the tail cover $\fb$ is {\em
  primitive}, otherwise we say that it is {\em new}.
\end{defn}

By \cite{Raynaud98}, Proposition 1.1.6, the converse of Condition (2)
holds automatically: if $\sigma< 1$ then there exists a unique
point $\bar{0}\in\Xb-\{\infb\}$ above which $\fb$ is tamely ramified.
The reason is this: if $\fb$ were \'etale over $\Xb-\{\infb\}$ and
$\sigma< 1$, then the Riemann--Hurwitz formula would produce a
negative number for the genus of $\Yb$.

It is clear from the results of \S \ref{simple} that the tail covers
obtained from the stable reduction of a three point Galois cover are
special (provided that the reduction is not exceptional). We can
therefore expect that a special tail cover $\fb:\Yb\to\Xb$ is `rigid'.
In fact, using the methods of \cite{BertinMezard00} one can show that
the curve $\Yb$ does not admit any nontrivial $G$-equivariant
deformations. What is more important for us is that all special
$G$-tail covers with a given ramification type $(m,h)$ are `locally
isomorphic', in the following sense.

\begin{defn}
  Let $(\fb_1:\Yb_1\to\Xb,\eta_1,H)$ and
  $(\fb_2:\Yb_2\to\Xb,\eta_2,H)$ be two pointed $G$-tail covers with
  the same subgroup $H\subset G$ and the same ramification type
  $(m,h)$. Let $\Xb_{\infb}$ denote the generic point of the completion
  of $\Xb$ at $\infb$ and set
  $\Yb_{i,\infb}:=\Yb_i\times_{\Xb}\Xb_{\infb}$. A {\em local
  isomorphism} between $\fb_1$ and $\fb_2$ is a $G$-equivariant
  isomorphism
  \[
         \varphi_{\infb}:\Yb_{1,\infb}\;\liso\;\Yb_{2,\infb}
  \]
  with the following properties.
  \begin{enumerate}
  \item
    The automorphism $\psi_{\infb}:\Xb_{\infb}\iso\Xb_{\infb}$ induced by
    $\varphi_{\infb}$ extends to a global automorphism
    $\psi:\Xb\iso\Xb$ (which necessarily fixes $\infb$).
  \item 
    We have $\vphi_{\infb}(\eta_1)=\eta_2$. Moreover, if $\fb_1$
    is ramified at $\bar{0}_1\in\Xb-\{\infb\}$ then $\fb_2$ is also
    ramified at some point $\bar{0}_2\in\Xb-\{\infb\}$, and the
    automorphism $\psi$ in (1) can be chosen in such a way that
    $\psi(\bar{0}_1)=\bar{0}_2$.
  \end{enumerate}
\end{defn}

\begin{lem} \label{Gtaillem}
  All special pointed $G$-tail covers with given ramification type
  $(m,h)$ and subgroup $H\subset G$ are locally isomorphic.
\end{lem}

\begin{proof}
This is a special case of the results of \cite{RachelWildly}, \S
4. For the convenience of the reader, we sketch a proof.

Let $(\fb:\Yb\to\Xb,\eta,H)$ be a special pointed $G$-tail cover, of
ramification type $(m,h)$. We will show that $\fb:\Yb\to\Xb$ is
uniquely determined, up to local isomorphism, by the datum
$(G,H,m,h)$. First of all, by the definition of local isomorphism it
suffices to look at tail covers $\fb$ which are totally ramified above
$\infb$. In other words, we may assume that $G=I\rtimes H$, where
$I\subset G$ is cyclic of order $p$ and $H$ is cyclic of order $m$.
Second, we shall assume that the tail cover $\fb$ is new. Hence we can
write $h=m+a$, with $0<a<m$. Note that $m|a(p-1)$ and that $a$ is
prime to $p$. We will indicate later how to modify the proof in the
case that $\fb$ is primitive.

Let $x$ be a global coordinate for $\Xb\cong\PP^1_k$ such that 
$\Xb_{\infb}=\Spec k((x^{-1}))$. By elementary Galois theory, the
localized cover $\Yb_{\infb}\to\Xb_{\infb}$ is given by two equations
\begin{equation} \label{Gtaillemeq1}
  z^m \;=\; x, \qquad 
  y^p-y \;=\; z^a\,(b_0 x + b_1 + b_2 x^{-1} + \ldots),
\end{equation}
with $b_i\in k$ and $b_0\not=0$. The action of $G$ on $\Yb_{\infb}$ is given
as follows. We may choose a generator $\alpha$ of $I\cong\ZZ/p$ such that
$\alpha^*z=z$ and $\alpha^*y=y+1$. Furthermore, for any $\beta\in H$ we have
$\beta^*z=\psi(\beta)\cdot z$ and $\beta^*y=\chi(\beta)\cdot y$, where
$\psi:H\to\bar{\FF}_p^\times$ is a character of order $m$ such that
$\chi=\psi^a$. See e.g.\ \cite{RachelWildly}, Lemma 1.4.1.

A local isomorphism amounts to a change of coordinates
\begin{equation} \label{Gtaillemeq2}
  x \;\mapsto\; c(x+d), \qquad 
  z \;\mapsto\; c^{1/m}z(1+dx^{-1})^{1/m}, \qquad
  y \;\mapsto\; y + g,
\end{equation}
where $c,d\in k$, $c\not=0$ and $g=z^a(e_1x^{-1}+e_2x^{-2}+\ldots)\in
k[[z^{-1}]]$. An easy computation shows that for a suitable choice of
the constants $c$ and $d$, and after transforming the equations
\eqref{Gtaillemeq1} by the change of coordinate \eqref{Gtaillemeq2},
we may assume that $b_0=1$ and $b_1=0$. It is also easy to see that we
may choose $e_1,e_2,\ldots$ in such a way that
\begin{equation} \label{Gtaillemeq3}
    g^p -g \;=\; -z^a\,(b_2x^{-1}+b_3x^{-2}+\ldots)
\end{equation}
(here we use that $m|a(p-1)$). With this choice, equations
\eqref{Gtaillemeq1} becomes 
\[
       z^m \;=\; x, \qquad y^p-y \;=\; z^h,
\]
which visibly depends only on the ramification type $(m,h)$. This
proves the lemma for a new tail cover $\fb$. The proof for a primitive
tail cover is similar. The main difference is that we are only allowed
to do a homothety $x\mapsto c x$ instead of a general affine
transformation $x\mapsto c(x+d)$. However, since $0<h<m$ the proof goes
through.  \end{proof}

\subsubsection{}  \label{Gdefdat3}

As before, $k$ denotes an algebraically closed field and $G$ a finite
group. Suppose we are given the following data.
\begin{itemize}
\item A subgroup $G_0\subset G$, of the form $G_0=I_0\rtimes H_0$,
  where $I_0\subset G_0$ is cyclic of order $p$ and $H_0\subset G_0$
  is of prime-to-$p$ order. Let $\chi_0:H_0\to \FF_p^\times$ denote
  the character defined by the rule
  $\beta\alpha\beta^{-1}=\alpha^{\chi_0(\beta)}$, where $\alpha\in
  I_0$, $\beta\in H_0$.
\item A normalized special deformation datum $(\Zb_0,\omega_0)$ over
  $\Xb_0:=\PP^1_k$ of type $(H_0,\chi_0)$. Let $(\tau_j)_{j\in B}$
  denote the set of critical points, $m_j,h_j$ and $\sigma_j$ the
  invariants attached to $(\Zb_0,\omega_0,\tau_j)$.
\item 
  For all $j\in B-B\wild$ a point $\xi_j\in\Zb_0$ above
  $\tau_j\in\Xb_0$.
\item 
  For all $j\in B-B\wild$ a subgroup $G_j\subset G$ and a
  connected pointed $G_j$-tail cover
  $(\fb_j:\Yb_j\to\Xb_j,\eta_j,H_j)$.
\end{itemize}
By abuse of notation, we shall write $(\Zb_0,\omega_0;\fb_j)$ to
denote all of the data introduced above.

\begin{defn} \label{Gdefdatdef}
  The datum $(\Zb_0,\omega_0;\fb_j)$ is called a {\em special
    $G$-deformation datum} if the following conditions hold.
  \begin{enumerate}
  \item
    The group $G$ is generated by its subgroups $G_0$ and $G_j$. 
  \item
    The tail cover $\fb_j:\Yb_j\to\Xb_j$ has ramification of type
    $(m_j,h_j)$. Moreover,
    \[
         I_{\eta_j} \;=\; I_0\rtimes H_j, \quad\text{\rm and}
           \quad H_j \;=\; H_0 \cap I_{\eta_j}.
    \]
  \end{enumerate}
\end{defn}

It is clear from the results of \S \ref{simple} that the stable
reduction of a three point cover with mild reduction gives rise to a
special $G$-deformation datum. To make this statement more precise, we
introduce the following terminology. 

\begin{defn} \label{stabGmapdef}
  A {\em stable $G$-map} over a scheme $S$ is a homomorphism $f:Y\to
  X$ of marked semistable curves over $S$ together with an $S$-linear
  action of $G$ on $Y$ such that $Y$ is stably marked and $f$ commutes
  with the $G$-action.  Two stable $G$-maps $f_i:Y_i\to X_i$ are said
  to be isomorphic if there exist isomorphisms $\vphi:Y_1\iso Y_2$ and
  $\psi:X_1\to X_2$ of marked curves such that
  $f_2\circ\vphi=\psi\circ f_1$ and such that $\vphi$ is
  $G$-equivariant.
\end{defn}

If $f:Y\to X$ is a $G$-cover over $\Kb$, as in \S \ref{stable}, then
we may consider the stable reduction $\fb:\Yb\to\Xb$ of $f$ as a
stable $G$-map over $k$.

Fix a special $G$-deformation datum $(\Zb_0,\omega_0;\Yb_j)$. We will
show that $(\Zb_0,\omega_0;\Yb_j)$ gives rise to a stable $G$-map
$\fb:\Yb\to\Xb$ over $k$. Define
\[
       \Xb \;:=\; (\;\Xb_0\;\; {\textstyle\coprod}\; \coprod_j
       \Xb_j\;)\,/_\sim.
\]
Here the equivalence relation `$\sim$' identifies the critical point
$\tau_j\in\Xb_0$ with the distinguished point $\infb_j\in\Xb_j$.  The
topological space $\Xb$ carries a unique structure of a semistable
$k$-curve such that $\Xb_0$ and $\Xb_j$ are the irreducible components
of $\Xb$. We define marked points $\xb_j$ on $\Xb$ (indexed by the set
$B_0=B\prim\cup B\wild$), as follows. For $j\in B\prim$ we let
$\xb_j\in\Xb_j$ be the unique tame branch point of the tail cover
$\fb_j$ (see the remark after Definition \ref{specialtaildef}). For
$j\in B\wild$ we set $\xb_j:=\tau_j$.

To define the curve $\Yb$, let $\Yb_0$ be the smooth projective curve
with function field $k(\Zb_0)^{1/p}$. The action of $H_0$ on $\Zb_0$
extends uniquely to an action on $\Yb_0$. We let $G_0$ act on $\Yb_0$
via the quotient map $G_0\to H_0$. Furthermore, we define
\[
    \Yb \;:=\; (\,\Ind_{G_0}^G(\Yb_0)\;\; {\textstyle\coprod}\;
                 \coprod_j\; \Ind_{G_j}^G(\Yb_j)\,)\,/_\sim.
\]
Here the equivalence relation `$\sim$' is defined as follows: for
$\eta'\in\Ind_{G_0}^G(\Yb_0)$ and $\eta''\in\Ind_{G_j}^G(\Yb_j)$ we have
$\eta'\sim \eta''$ if and only if there exists an element $\alpha\in
G$ such that $\alpha(\eta')\in\Yb_0$ lies above $\xi_j\in\Zb_0$ and
$\alpha(\eta'')=\eta_j\in\Yb_j$. (This is an equivalence relation
because of Condition (2) of Definition \ref{Gdefdatdef}.) Again,
$\Yb$ carries a canonical structure of a semistable curve over $k$. It
follows from Condition (1) of Definition \ref{Gdefdatdef} that $\Yb$
is connected. Let $\fb:\Yb\to\Xb$ be the unique map which commutes
with the natural $G$-action on $\Yb$ and whose restriction to $\Yb_0$
(resp.\ to $\Yb_j$) is equal to $\fb_0$ (resp.\ to $\fb_j$). The
marked points of $\Yb$ are, by definition, those which lie over a
marked point of $\Xb$. Thus, a point on $\Yb$ is marked either if it
is the $G$-translate of a ramification point of a tail cover $\fb_j$
or if it is a $G$-translate of a point on $\Yb_0$ which lies above a
wild critical point. One checks that $\fb:\Yb\to\Xb$ is a stable
$G$-map. We call $\Xb_0=\PP^1_k$ the original component of $\fb$.
 
\begin{defn} \label{specialGmapdef}
  A stable $G$-map $\fb:\Yb\to\Xb$ over $k$ is called {\em special} if
  there exists a special $G$-deformation datum $(\Zb_0,\omega_0;\Yb_j)$
  such that $\fb$ is isomorphic to the stable $G$-map constructed
  above.
\end{defn} 

We can summarize the results of \S \ref{stablered} and \S \ref{simple}
as follows.

\begin{prop} \label{Gdefdatprop}
  Let $f:Y\to X=\PP^1_{\Kb}$ be a $G$-Galois cover branched at $0$,
  $1$ and $\infb$, with mild reduction. Then the stable reduction of
  $f$ is either exceptional or special.
\end{prop}

In \S \ref{lifts} we will prove the `converse' of this proposition:
every special $G$-map arises as the stable reduction of a three point
cover $f:Y\to X=\PP^1_{\Kb}$ with mild reduction.

\subsubsection{} \label{Gdefdat4}

Let us fix a special $G$-deformation datum $(\Zb_0,\omega_0;\Yb_j)$
and let $\fb:\Yb\to\Xb$ be the associated special $G$-map. Let us
denote by $\Aut_G^0(\fb)$ the group of automorphisms of the
special $G$-map $\fb$ which induce the identity on the original
component $\Xb_0$. If $\fb$ is the stable reduction of a three point
$G$-cover $f:Y\to X$ then the image of the inner monodromy action
$\kappab\inn$ associated to $f$ is contained in
$\Aut_G^0(\fb)/C_G$. We are interested in a description of
$\Aut_G^0(\fb)$ in terms of the special deformation datum
$(\Zb_0,\omega_0;\Yb_j)$. 

\begin{lem} \label{autflem}
  There is a natural bijection between the group $\Aut_G^0(\fb)$ and
  the set of tuples $(\gamma_0;\gamma_j)$, where $\gamma_0\in G$ and
  $\gamma_j:\Yb_j\iso\Yb_j$ is an (outer) automorphism of the pointed
  $G_j$-tail cover $(\fb_j,\eta_j,H_j)$ (see Definition
  \ref{tailcoverdef}) such that the following holds.
  \begin{enumerate}
  \item 
    The element $\gamma_0\in G$ normalizes
    $I_0$ and centralizes $H_0$.
  \item
    The equation
    \[
       \gamma_0^{-1}\circ\alpha\circ\gamma_0 \;=\;
          \gamma_j\circ\alpha\circ\gamma_j^{-1}
    \]
    holds for all $\alpha\in G_j$ and $j\not\in B\wild$.
  \end{enumerate}
\end{lem}

\begin{proof} Let $\gamma:\Yb\iso\Yb$ be an element of $\Aut_G^0(\fb)$. It is
clear that we can find an element $\gamma_0\in G$ such that
$\gamma_0^{-1}\circ\gamma$ induces the identity on $\Yb_0$. It is also
clear that $\gamma_0$ normalizes $G_0$ and that conjugation by
$\gamma_0$ induces the identity on $G_0/I_0\cong H_0$. Since the order
of $H_0$ is prime to $p$, the Schur--Zassenhaus Theorem implies that the
subgroup $\gamma_0H_0\gamma_0^{-1}\subset G_0$ is conjugate to $H_0$
by an element of $G_0$. Using the fact that the action of $H_0$ on
$I_0$ is irreducible, one shows that there exists a unique element
$\alpha\in I_0$ such that
\[
      \gamma_0\circ\beta\circ\gamma_0^{-1} \;=\; 
         \alpha^{-1}\circ\beta\circ\alpha
\]
for all $\beta\in H_0$. After replacing $\gamma_0$ by
$\alpha\circ\gamma_0$, we may thus assume that $\gamma_0$ centralizes
$H_0$. For all $j\not\in B\wild$, set
$\gamma_j:=\gamma_0^{-1}\circ\gamma|_{\Yb_j}$. One easily checks that
$\gamma_j$ is an (outer) automorphism of the pointed $G_j$-tail cover
$(\fb_j,\eta_j,H_j)$ and that the tuple $(\gamma_0;\gamma_j)$
satisfies Condition (1) and (2).

Conversely, given a tuple $(\gamma_0;\gamma_j)$ such that (1) and (2)
hold, we define an automorphism $\gamma:\Yb\iso\Yb$ of
$\fb$, as follows. If $y$ is a point on $\Yb_0$ and $\alpha$ is an
element of $G$, we set
\[
      \gamma(\alpha(y)) \;:=\; \alpha\circ\gamma_0(y).
\]
This is well defined because of Condition (1).  If $y$ is a point on
$\Yb_j$ and $\alpha\in G$, we set
\[
     \gamma(\alpha(y)) \;:=\; \alpha\circ\gamma_0\circ\gamma_j(y).
\]
This is well defined because $\gamma\circ\gamma_j$ commutes with
$G_j$, by Condition (2). Moreover, the two definitions respect the
equivalence relation `$\sim$' used in the definition of $\Yb$ in \S
\ref{Gdefdat3}. Hence we obtain an automorphism $\gamma:\Yb\iso\Yb$,
which is easily seen to be an element of $\Aut_G^0(\fb)$. This proves
the lemma.  \end{proof}

The association $\gamma\mapsto \gamma_0$ is a group homomorphism
$\Aut_G^0(\fb)\to G$. Let us denote by $n$ the order of the image of
this homomorphism in $G/C_G$. Note that $n|(p-1)$. Let us denote by
$\Aut_G(\fb_j)$ the group of inner automorphisms of the pointed
$G_j$-tail covers $(\fb_j,\eta_j,H_j)$. The lemma shows that
\[
      |\Aut_G^0(\fb)| \;=\; n\cdot|C_G|\cdot\prod_{j\notin B\wild}\,
                 |\Aut_G(\fb_j)|.
\]

\begin{rem} \label{autfrem}
  Suppose that the character $\chi:H_0\to\FF_p^\times$ is
  injective. (In particular, $H_0$ is cyclic). Then the integer $n$
  defined above is a multiple of 
  \[
        n' \;:=\; {\rm gcd}_j(m_j).
  \]
  In fact, let $\gamma$ be an element of $H_0$ which is contained in
  the stabilizer $H_j$ of $\xi_j$, for all $j\in B\prim\cup B\new$.
  Then the tuple $(\gamma_0;\gamma_j)$ with
  $\gamma_j:=\gamma_0^{-1}|_{\Yb_j}$ represents an element of
  $\Aut_G^0(\fb)$ which does not lie in $C_G$ and which acts trivially
  on the vertical components of $\Yb$. 
\end{rem}


\section{The auxiliary cover} \label{aux}

Let $(\Zb,\omega)$ be a special deformation datum of type $(H,\chi)$.
By a {\em lift} of $(\Zb,\omega)$ we mean a Galois cover $f:Y\to
X=\PP^1$ with Galois group $G:=\ZZ/p\rtimes_\chi H$ such that the
deformation datum occuring over the original component of the stable
reduction of $f$ is equal to $(\Zb,\omega)$.  As we will see below,
such a lift is uniquely determined by the deformation datum
$(\Zb,\omega)$ and the set of branch points $(x_j)$ of $f$ (which lift
the critical points $(\tau_j)$ of $(\Zb,\omega)$).

Following \cite{special}, we say that the lift $f:Y\to X$ is {\em
special} if the stable reduction of $f$ is as simple as possible,
i.e.\ if the curve $\Xb$ is a comb. For instance, if $f:Y\to X$ is the
{\em auxiliary cover} associated to a three point cover with mild
reduction, then $f$ is special by Theorem \ref{thm1} (note that $f$ is
in general not a three point cover itself). 

In \cite{special} we showed that a lift $f$ of $(\Zb,\omega)$ is
special if each branch point of $f$ is `sufficiently closed' to a
$K_0$-rational point. The main result of this section -- Theorem
\ref{stabauxthm} -- says that the converse is true as well. This
theorem will play a crucial role in \S \ref{main}.

The proof of Theorem \ref{stabauxthm} relies on the results of
\cite{cotang}, which concern the deformation theory of a certain curve
with group scheme action associated to $(\Zb,\omega)$. The connection
is as follows. If $f:Y\to X$ is a lift of the deformation datum
$(\Zb,\omega)$ then there exist a certain natural $R$-model
$f_R:Y_R\to X_R=\PP^1_R$ (which is {\em not} the stable model) such
that the $G$-action on $Y$ extends to an action of the group scheme
$\G:=\bmu_p\rtimes H$ on $Y_R$. The special fiber $\Yb$ of $Y_R$ is a
singular curve with an action of $\G$; it depends on $(\Zb,\omega)$
but not on the lift $f:Y\to X$. Therefore, a lift of $(\Zb,\omega)$ is
the same thing as an equivariant deformation of $\Yb$.

\subsection{The $\G$-cover associated to a multiplicative 
            deformation datum} \label{defo}

\subsubsection{} \label{defo1}

Throughout this section, we fix an odd prime $p$ and an
algebraically closed field $k$ of characteristic $p$. We denote by
$W(k)$ the ring of Witt vectors over $k$. 

Let $\Xb$ be a smooth projective curve over $k$ and let $(\Zb,\omega)$
be a multiplicative deformation datum over $\Xb$ of type $(H,\chi)$,
see Definition \ref{defdatdef}. In this section we show that
$(\Zb,\omega)$ corresponds naturally to a (not necessarily smooth)
curve $\Yb$ over $k$ together with an action of a certain group scheme
$\G$ such that $\Xb=\Yb/\G$.

We define $\G$ as a group scheme over $W(k)$, as follows. If $R$ is a
$W(k)$-algebra, an element of $\G(R)$ is given by a tuple
$(\zeta,\beta)$, where $\zeta\in R$ satisfies $\zeta^p=1$ and
$\beta\in H$. Multiplication is defined as follows:
\[
    (\zeta_1,\beta_1)\cdot(\zeta_2,\beta_2) \;:=\; 
       (\zeta_1\zeta_2^{\chi(\beta_1)},\beta_1\beta_2).
\]
In other words: $\G$ is the semidirect product $\bmu_p\rtimes H$
determined by the action of $H$ on $\bmu_p$ via $\chi$. If $R$ is a
$W(k)$-algebra and $Y$ a scheme over $R$, we say that $\G$ acts on $Y$
if we are given an $R$-linear action of $\G\otimes R$ on $Y$.

Let $(\tau_j)_{j\in B}$ be the set of critical points for the
deformation datum $(\Zb,\omega)$. For $j\in B$, let $m_j,h_j$ be the
invariants attached to $(\Zb,\omega,\tau_j)$, see \S \ref{defdat}. A
critical point $\tau_j$ is called a {\em tame} (resp.\ {\em wild})
{\em branch point} if $m_j>1$ (resp.\ if $h_j=0$). We denote by
$B\tame $ and $B\wild$ the corresponding subsets of $B$ and set
$B\branch:=B\tame\cup B\wild$. In general we have $B\tame\cap
B\wild\not=\emptyset$ and $B\branch\not=B$.

Let $\Ub\subset\Xb$ be the complement of the wild branch points and
let $\Vb\subset\Zb$ be its inverse image. Then $\omega$ is a regular
and logarithmic differential form on $\Vb$. Let $\Wb\to\Vb$ be the
corresponding $\bmu_p$-torsor, see e.g.\ \cite{MilneEC}, p. 129.  The
definition of $\Wb$ is roughly the following. Locally on $\Vb$ we can
write $\omega=\diff u/u$ for a unit $u$. Then $y^p=u$ is a local
equation for the $\bmu_p$-torsor $\Wb\to\Vb$.

\begin{prop} \label{defo1prop}
  Let $\Yb$ be the projective compactification of $\Wb$ such
  that the complement $\Yb-\Wb$ is contained in the smooth locus of
  $\Yb$. The action of $\bmu_p$ on $\Wb$ extends uniquely to an action
  of $\G$ on $\Yb$ such that $\Zb=\Yb/\bmu_p$ and $\Xb=\Yb/\G$.
\end{prop}

\begin{proof} Let $\lambda:\bmu_p\times\Wb\to\Wb$ denote the action of
$\bmu_p$ on $\Wb$. We have to show that $\lambda$ extends to all of
$\Yb$. Let $\eta\in\Yb-\Wb$ be a point in the complement of $\Wb$ and
let $\xi\in\Zb$ be its image. By construction, the point $\xi$ lies
above a wild critical point $\tau_j$, so $\omega$ has a simple pole at
$\xi$. Write $\omega=a\,\diff u/u$, where $u$ is a regular function in
a neighborhood of $\xi$ with a {\em simple} zero at $\xi$ and
$a\in\FF_p^\times$ is the residue of $\omega$ at $\xi$. By definition
of $\Yb$ there exists a rational function $y$ on $\Yb$ such that
$y^p=u$ and $\lambda^*y=\zeta^b\otimes y$, for some
$b\in\FF_p^\times$. It follows also from the definition of $\Yb$ that
$\eta$ is a smooth point of $\Yb$. Therefore, the integral equation
$y^p=u$ shows that $y\in\OO_{\Yb,\eta}$. Moreover, since $u$ has a
simple zero at $\xi$, we have $\OO_{\Yb,\eta}=\OO_{\Zb,\xi}[y\mid
y^p=u]$.  This shows that the action $\lambda$ extends to $\Yb$.
It is obviously that this extension is unique.

It remains to show that $\lambda$ extends to an action of $\G$ on
$\Yb$ which induces the canonical action of $H$ on $\Zb$. But this is
clear from the rule $\beta^*\omega=\chi(\beta)\cdot\omega$ and the
definition of $\Yb$.  \end{proof}

The map $\Yb\to\Xb$ together with the action of $\G$ on $\Yb$ will be
called the $\G$-cover associated to $(\Zb,\omega)$. Note that the
curve $\Yb$ is integral and has isolated, `cusp-like' singularities
precisely above the critical points which are not wild branch
points, see \cite{Raynaudfest}.

\subsubsection{} \label{defo2}

Let $\C_k$ be the category of Noetherian Artinian local $W(k)$-algebras
with residue field $k$. An {\em equivariant deformation} of $\Yb$ over
$R\in\C_k$ is a deformation $Y$ of the curve $\Yb$ over $R$ together
with an action of $\G$ on $Y$ such that the natural isomorphism
$\Yb\cong Y\otimes_R k$ is $\G$-equivariant. Let $\Def(\Yb,\G)$ denote
the functor which associates to a ring $R\in\C_k$ the set of
isomorphism classes of equivariant deformations of $\Yb$ over $R$.

Given an equivariant deformation $Y_R$ of $\Yb$, we denote by
$Z_R:=Y_R/\bmu_p$ and $X_R:=Y_R/\G$ the quotient schemes by $\bmu_p$
and $\G$, respectively. The $R$-curve $X_R$ (resp.\ $Z_R$) is a smooth
deformation of the curve $\Xb$ (resp.\ $\Zb$).  The association
$Y_R\mapsto X_R$ gives rise to a morphism of deformation functors
$\Def(\Yb,\G)\to\Def(\Xb)$. Here $\Def(\Xb)$ denotes the functor which
classifies deformations of $\Xb$. Let us write $\Def(\Xb;\tau_j\mid
j\in B\branch)$ for the functor which classifies deformations of the
pointed curve $(\Xb;\tau_j\mid j\in B\branch)$. We claim that the
morphism $\Def(\Yb,\G)\to\Def(\Xb)$ factors in a natural way over the
forgetful morphism $\Def(\Xb;\tau_j\mid j\in B\branch)\to\Def(\Xb)$,
thus inducing a morphism
\begin{equation} \label{functorhom}
   \Def(\Yb,\G) \;\To\; \Def(\Xb;\tau_j\mid j\in B\branch).
\end{equation}

Let $Y_R$ be an equivariant deformation of $\Yb$ and set
$X_R:=Y_R/\G$, $Z_R:=Y_R/\bmu_p$.  To prove the claim we have to
endow, in a functorial way, the curve $X_R$ with sections
$\tau_{j,R}:\Spec R\to X_R$ lifting the branch points $\tau_j$, for
all $j\in B\branch$. It is clear that the natural map $Z_R\to X_R$ is
a tame $H$-Galois cover lifting $\Zb\to\Xb$. By the theory of tame
ramification, the branch divisor of $Z_R\to X_R$ consists of sections
$\tau_{j,R}:\Spec R\to X_R$ lifting the tame branch points
$\tau_j\in\Xb$, for all $j\in B\tame$. Now let $j\in B\wild$ and
choose a point $\xi\in\Zb$ above $\tau_j$. We can write $\omega=\diff
u/u$ with $\ord_\xi(u)=1$. By construction, the restriction of
$\Yb\to\Zb$ to a neighborhood of $\xi$ is a $\bmu_p$-cover given by
the Kummer equation $y^p=u$. Hence the deformation $Y_R\to Z_R$ is,
locally around $\xi$, still a $\bmu_p$-cover given by a Kummer
equation $y^p=u_R$, where $u_R$ lifts $u$. The equation $u_R=0$
defines a section $\xi_R:\Spec R\to Z_R$ lifting $\xi$. We define
$\tau_{j,R}:\Spec R\to X_R$ as the image of $\xi_R$. Using the
$H$-action, it is easy to see that our two definitions of $\tau_{j,R}$
agree if $j\in B\tame\cap B\wild$. This proves the claim.

We say that $Y_R\to X_R$ is a deformation of the $\G$-cover
$\Yb\to\Xb$, and we call the tuple $(\tau_{j,R})_{j\in B\branch}$ the
{\em branch locus} of $Y_R\to X_R$. It is clear that a deformation of
the $\G$-cover $\Yb\to\Xb$ is the same thing as an equivariant
deformation of $\Yb$.

\subsection{The case of a special deformation datum} 
       \label{specialcase}

\subsubsection{} \label{specialcase1}

We continue with the notation of \S \ref{defo}. In addition, we assume
from now on that the deformation datum $(\Zb,\omega)$ is {\em
  special}, see Definition \ref{specialdefdatdef}. We use the notation
$B\prim$, $B\new$ and $B_0$ as in \S \ref{Gdefdat1}. It follows
immediately from the definitions that we have $B\tame=B\prim\cup
B\new$ and $B\wild=B_0-B\prim$. In particular, we have $B\tame\cap
B\wild=\emptyset$ and $B=B\branch$. Furthermore, we assume that the
special deformation datum $(\Zb,\omega)$ is normalized, i.e.\ we have
$\Xb=\PP^1_k$ and $\{\tau_j\mid j\in B_0\}=\{0,1,\infty\}$.

Let $Y_R\to X_R$ be a deformation of the $\G$-cover $\Yb\to \Xb$ with
branch locus $(\tau_{j,R})_{j\in B}$. There exists a unique
isomorphism $X_R\cong\PP^1_R$ which reduces to the identity on
$\Xb=\PP^1_k$ and such that $\{\,\tau_{j,R}\mid j\in
B_0\,\}=\{0,1,\infty\}$. Hence we may and will identify $X_R$ with
$\PP^1_R$ and the sections $\tau_{j,R}$, for $j\in B\new$, with the
corresponding elements of $R$. A $B$-tuple $(\tau_{j,R})$ of
$R$-points of $X_R=\PP^1_R$ such that $\tau_{j,R}$ lifts $\tau_j$ for
all $j\in B$ and $\{\,\tau_{j,R}\mid j\in B_0\,\}=\{0,1,\infty\}$ is
called a {\em normalized $B$-tuple over $R$}.

\begin{prop} \label{specialcase1prop}
  The morphism of deformation functors 
  \[
      \Def(\Yb,\G) \;\To\; \Def(\Xb;\tau_j)
  \]
  is an isomorphism.  In other words: given $R\in\C_k$ and a
  normalized $B$-tuple $(\tau_{j,R})$ over $R$, there exists a
  deformation $Y_R\to X_R=\PP^1_R$ of the $\G$-cover
  $\Yb\to\Xb=\PP^1_k$ with branch locus $(\tau_{j,R})$. Moreover,
  $Y_R\to X_R$ is unique up to unique isomorphism.
\end{prop}

\begin{proof}
This is a special case of \cite{cotang}, Theorem 5.7.
\end{proof}

For $j\in B\new$, let $T_j$ be an indeterminate; set
$\Rt:=W(k)[[\,T_j\mid j\in B\new\,]]$. For $j\in B\new$, let
$[\tau_j]\in W(k)$ denote the Teichm\"uller lift of $\tau_j$ (in fact,
any lift of $\tau_j$ to an element of $W(k)$ would do as well). Let
$(\tilde{\tau}_j)$ be the normalized $B$-tuple over $\Rt$ which has
$\tilde{\tau}_j=[\tau_j]+T_j$ for $j\in B\new$. Let $\X$ denote the
formal completion along the special fiber of the projective line over
$\Rt$. It is clear that for any $R\in\C_k$ there is a one-to-one
correspondence between local $W(k)$-algebra homomorphism $\Rt\to R$
and normalized $B$-tuples over $R$. Therefore, $(\X;\tilde{\tau}_j)$
is the universal formal deformation of the pointed curve
$(\Xb;\tau_j)$. Hence Proposition \ref{specialcase1prop} implies:

\begin{cor} \label{specialcase1cor} 
  Let $\Y\to\X$ be the (unique) formal deformation of the $\G$-cover
  $\Yb\to\Xb$ with branch locus $(\tilde{\tau}_j)$. Then $\Y$ is a
  universal formal equivariant deformation of $\Yb$.
\end{cor}

\subsubsection{} \label{specialcase2}

For $j\in B$, we denote by $\Xd_j$ the completion of $\Xb$ at $\tau_j$
and set $\Yd_j:=\Yb\times_{\Xb}\Xd_j$. The action of $\G$ on $\Yb$
induces an action on $\Yd_j$ such that $\Xd_j=\Yd_j/\G$. Given a
deformation $Y_R\to X_R$ of the $\G$-cover $\Yb\to\Xb$, we denote by
$\Xd_{j,R}$ the completion of $X_R$ at $\tau_j$ and set
$\Yd_{j,R}:=\Xd_j\times_{X_R}\Xd_{j,R}$. Then $\Yd_{j,R}$ is a
$\G$-equivariant deformation of $\Yd_j$. The association
$Y_R\mapsto(\Yd_{j,R})_{j\in B}$ defines a morphism of deformation
functors
\begin{equation} \label{locglobmorph}
   \Phi:\,\Def(\Yb,\G) \;\To\; \prod_{j\in B}\;\Def(\Yd_j,\G),
\end{equation}
Here $\Def(\Yd_j,\G)$ denotes the deformation functor which
classifies $\G$-equivariant formal deformations of $\Yd_j$. Following
\cite{BertinMezard00}, we call $\Phi$ the {\em local-global morphism}.

\begin{prop} \label{specialcase2prop}
  \begin{enumerate}
  \item
    The functor $\Def(\Yd_j,\G)$ admits a miniversal deformation
    over the ring
    \[
       \Rt_j \;:=\quad\begin{cases}
            \quad W(k),\qquad & \text{\rm for $j\in B_0$,} \\
            \quad W(k)[[t_j]],\qquad & \text{\rm for $j\in B\new$.}
           \end{cases}
    \]
    Here $t_j$ is an indeterminate. 
  \item
    The local-global morphism $\Phi$ is an isomorphism. Therefore,
    \[
       \Rt \;\cong\; \widehat{\otimes}_{W(k)}\, \Rt_j \;=\;
         W(k)[[\,t_j\mid j\in B\new\,]].
    \]
  \item
    For every $j\in B\new$ there exists a unit $w_j\in\Rt$ such that
    \[
            t_j \;\equiv\; w_jT_j \pmod{p}.
    \]
  \end{enumerate}
\end{prop}

\begin{proof}
Assertion (1) and (2) are a special case of \cite{cotang}, Theorem
5.11. Assertion (3) corresponds to \cite{cotang}, Proposition 5.14. 
\end{proof}

\subsubsection{} \label{specialcase3}

Let $\Y\to\X$ be the universal deformation over $\Rt$ of the
$\G$-cover $\Yb\to\Xb$, as defined in \S \ref{specialcase1}.  Fix an
index $j\in B\new$. We denote by $\hat{\X}_j$ the completion of $\X$
at $\tau_j$ and set $\hat{\Y}_j:=\Y\times_{\X}\hat{\X}_j$. Being an
equivariant deformation of $\Yd_j$ over $\Rt$, $\hat{\Y}_j$ is induced
from the miniversal deformation of Proposition \ref{specialcase2prop}
via a homomorphism of $W(k)$-algebras $\varphi_j:\Rt_j\to\Rt$. The
homomorphism $\varphi_j$ is unique up to an automorphism of
$W(k)$-algebra $\Rt_j$. Since $\varphi_j$ is injective by Proposition
\ref{specialcase2prop} (2), we may identify the ring $\Rt_j$ with the
image of $\varphi_j$. By Proposition \ref{specialcase2prop} (1) the
ring $\Rt_j$ is a formal power series over $W(k)$ in one variable
$t_j$. It is clear that Assertion (3) of Proposition
\ref{specialcase2prop} does not depend on the choice of the parameter
$t_j$ as long as $\Rt_j\cong W(k)[[t_j]]$.

The quotient $\Z:=\Y/\bmu_p$ is a formal $H$-equivariant deformation
of $\Zb$. Let us choose a point $\xi\in\Zb$ above $\tau_j$ and denote
by $\hat{\Z}_j$ the completion of $\Z$ at $\xi_j$. Since $\Z$ is
formally smooth over $\Rt$ we have $\hat{\Z}_j=\Spec\Rt[[z]]$ for some
parameter $z$.  We may assume
\begin{equation} \label{zeq}
   \beta^*z \;\equiv\; \chi_j(\beta)\cdot z \pmod{p},
\end{equation}
for all $\beta\in H_j$ and some character $\chi_j:H_j\to k^\times$.
Actually, we have $\chi_j^{h_j}=\chi|_{H_j}$. (Of course, we could
have choosen $z$ in such a way that the above congruence would become
an equality in $\Rt[[z]]$, but this would give a stronger restriction
on the choice of $z$ than we are are willing to make.)

The natural map $\Y\to\Z$ is a $\bmu_p$-torsor in a neighborhood of
$\xi$. Therefore, the pullback to the completion $\hat{\Z}_j$ is
defined by a Kummer equation
\begin{equation} \label{Kummereq1}
      y^p \;=\; u \;=\; \sum_{\mu\geq 0} \, c_\mu\,z^\mu,
\end{equation}
with $c_\mu\in \Rt$ and $c_0\in\Rt^\times$. The unit
$u\in\Rt[[z]]^\times$ is unique up to a $p$th power. As we have
already remarked above, the $\bmu_p$-torsor $\hat{\Y}_j\to\hat{\Z}_j$
is the pullback of a $\bmu_p$-torsor defined over $\Rt_j$. In other
words, we can choose the parameter $z$ and the Kummer equation
\eqref{Kummereq1} in such a way that all coefficients $c_\mu$ lie in
$\Rt_j\subset\Rt$. In the following, we will assume that
$c_\mu\in\Rt_j$ and we pretend that $\hat{\Y}_j\to\hat{\Z}_j$ is
actually defined over $\Rt_j$.

\begin{lem} \label{specialcase3lem}
  Let $a_j:=h_j-m_j$. 
  We may choose the parameter $z$, the Kummer equation
  \eqref{Kummereq1} and a parameter $t_j$ for $\Rt_j$ such that 
    \[
       u \;\equiv\; 1 \;+\; t_j\,z^{a_j} \;+\; z^{h_j} 
         \mod{(\,p,\,t_j\,z^{a_j+1},\,z^{h_j+1}\,)}.
  \]
\end{lem}

\begin{proof} Given an element $c\in\Rt_j$, we write $c'$ (resp. $\bar{c}$)
for the image of $c$ in $R':=\Rt_j/(p)$ (resp.\ in $k=\Rt_j/(p,t_j)$).
We claim that $c_0'\in({R'}^\times)^p$. To prove this claim, let
$\beta$ be a generator of $H_j\cong\ZZ/m_j$ and $0<a<p$ an integer
such that $a\equiv\chi(\beta)\mod{p}$. The subgroup scheme
$\G_j:=\bmu_p\cdot H_j\subset\G$ acts on $\hat{\Y}_j$, extending the
action of $H_j$ on $\hat{\Z}_j$ given by \eqref{zeq}. Therefore, there
exists a unit $w\in\Rt_j[[z]]^\times$ such that $\beta^*y=wy^a$.
Taking $p$th powers on both sides we obtain $\beta^*u=w^pu^a$. Using
\eqref{zeq} shows that $c_0^{1-a}$ is congruent to a $p$th power
modulo $p$. Since $\sigma_j=h_j/m_j$ is not an integer, we have
$a\not\equiv 1\mod{p}$. The claim follows. In particular, after
multiplying $u$ with a $p$th power in $\Rt_j[[z]]$ we may assume that
$c_0'=1$. We may also assume that $c_\mu'=0$ whenever $p|\mu$.

Let $\omega':=\diff u'/u'$.  This is a logarithmic differential which
lifts $\omega$ from $\Zd_j$ to $\hat{\Z}_j\otimes R'$; it classifies
the $\bmu_p$-torsor $\hat{\Y}_j\otimes R'\to\hat{\Z}_j\otimes R'$ up
to isomorphism. Since the action of $\bmu_p$ on $\hat{\Y}_j$ extends
to an action of $\G_j$, we have
$\beta^*\omega'=\chi(\beta)\cdot\omega'$ for all $\beta\in H_j$. This
and \eqref{zeq} shows that
\[
     \omega' \;=\; 
       \sum_{\mu\geq 0}\,b_\mu'\,z^{\mu m_j+a_j-1}\,\diff z
\]
with $b_\mu'\in R'$ (note that $0<a_j<m_j$ and $\chi_j^{a_j}=
\chi|_{H_j}$). Since $\omega$ has a zero of exact order $h_j-1$ at
$\xi_j$, we have $\bar{b}_0=0$ and $\bar{b}_1\not=0$. After a short
computation, we conclude that $c_\mu'=0$ for $0<\mu<a_j$ and
$\bar{c}_\mu=0$ for $a_j<\mu<h_j$. 

So far, we have shown that
\[
     u \;\equiv\; 1 \;+\; ct_j\,z^{a_j} \;+\; z^{h_j} 
         \quad\mod{(\,p,\,t_j\,z^{a_j+1},\,z^{h_j+1}\,)}
\]
for some element $c\in\Rt_j$.  It remains to prove that
$\bar{c}\not=0$. For then we can replace $t_j$ by $c^{-1}t_j$ and
obtain the desired congruence. Let $k[\epsilon]$ be the ring of dual
numbers over $k$ and $\kappa:R_j\to k[\epsilon]$ 
the unique $W(k)$-algebra homomorphism with $\kappa(t_j)=\epsilon$.
We can write
\[
   \kappa(u) \;=\; \bar{u}+\epsilon\cdot\bar{v} \;=\;
      (1 + z^{h_j} + \sum_{\mu\geq h_j}\,\bar{c}_\mu\,z^\mu)  
         + \epsilon\cdot(\,\sum_{\mu\geq a_j}\,\bar{d}_\mu\,z^\mu),
\]
with $\bar{c}_\mu,\bar{d}_\mu\in k$. We have to show that
$\bar{d}_{a_j}\not=0$.

By the proof of \cite{cotang}, Lemma 5.12, there exists an
$H_j$-invariant derivation of $k[[z]][z^{-1}]$
\[
   \theta_j \;=\; (\bar{e}_0+\bar{e}_1\,x+\ldots)\,\vf{x} \;=\;
     \frac{z^{1-m_j}}{m_j}\,(\,\bar{e}_0+\bar{e}_1\,z^{m_j}
                +\ldots)\,\vf{z}
\]
(which is the same thing as a derivation of $k[[x]]$, with
$x:=z^{m_j}$) such that the following holds:
\begin{itemize}
\item[(a)]
  $\bar{v}=\theta_j(\bar{u})$,
\item[(b)]
  we have $\bar{e}_0=0$ if and only if the pullback of $\hat{\Y}_j$
  via the homomorphism $\kappa:\Rt_j\to k[\epsilon]$ is the
  trivial deformation of $\Yd_j$. 
\end{itemize}
In fact, the class of $\theta_j$ modulo $k[[x]]\,x\,\vf{x}$
represents the isomorphism class of the equivariant deformation
$\hat{\Y}_j\otimes k[\epsilon]$, via a certain isomorphism
\[
    \EXt_{\G}^1(\Omega_{\Yb},\OO_{\Yb})_{\tau_j} \;\cong\;
       \T_{\Xb,\tau_j}/x\cdot\T_{\Xb,\tau_j},
\]
see \cite{cotang}, \S 5.3. But since $\hat{\Y}_j$ is a miniversal
deformation of $\Yd_j$ and $\kappa$ is not constant,
$\hat{\Y}_j\otimes k[\epsilon]$ must be nontrivial. Therefore, 
\[
      \bar{d}_{a_j} \;=\; \frac{h_j\,\bar{e}_0}{m_j} \;\not=\; 0.
\]
This proves the lemma.  \end{proof}

\begin{cor}
  We may choose the parameter $z$ and the Kummer equation
  \eqref{Kummereq1} in such a way that
  \[
       u \;\equiv\; 1 \;+\; w_jT_j\,z^{a_j} \;+\; z^{h_j} 
         \mod{(\,p^2,\,T_j\,z^{a_j+1},\,z^{h_j+1}\,)}.
  \]
  for some unit $w_j\in\Rt^\times$. 
\end{cor}

\begin{proof}
It  follows from Proposition \ref{specialcase2prop} (3) and Lemma
\ref{specialcase3lem} that we may assume
\[
    u \;\equiv\; 1 \;+\; w_jT_j\,z^{a_j} \;+\; z^{h_j} + p\cdot v
         \mod{(\,p^2,\,T_j\,z^{a_j+1},\,z^{h_j+1}\,)}
\]
for a unit $w_j\in\Rt$ and some $v\in\Rt[[z]]$. Replacing the
parameter $z$ by $z':=z-pv(h_j+a_jw_jT_j)^{-1}$, we obtain the stronger
congruence, modulo $p^2$.
\end{proof}

\subsection{The stable reduction of the auxiliary cover} \label{stabaux}

\subsubsection{} \label{stabaux1}

Let $k$ be an algebraically closed field of characteristic $p>0$ and
let $(\Zb,\omega)$ be a normalized special deformation datum over
$\Xb:=\PP^1_k$ of type $(H,\chi)$. Let $\Yb\to\Xb$ be the $\G$-cover
associated to $(\Zb,\omega)$, as defined in \S \ref{defo1}. Let $K_0$
denote the fraction field of $W(k)$ and let $K/K_0$ be a finite
extension. We assume that $K$ contains the $p$th roots of unity. The
integral closure $R$ of $W(k)$ inside $K$ is a complete mixed
characteristic discrete valuation ring with residue field $k$. 

Let $(\tau_{j,R})$ be a normalized $B$-tuple over $R$, see \S
\ref{specialcase1}. By Corollary \ref{specialcase1cor} there exists a
unique formal deformation $Y_R\to X_R=\PP^1_R$ of the $\G$-cover
$\Yb\to\Xb$ with branch locus $(\tau_{j,R})$. By Grothendieck's
Existence Theorem, this formal deformation is effective, i.e.\ we may
regard $Y_R\to X_R$ as a finite morphism of schemes over $R$. Write
$Y_K\to X_K=\PP^1_K$ for the generic fiber of $Y_R\to X_R$ and
$x_j:=\tau_{j,K}$ for the $K$-rational point on $X_K$ corresponding to
the section $\tau_{j,R}$, for $j\in B$. Let
\[
   G \;:=\; \G(K) \;\cong\; \ZZ/p\rtimes_{\chi}H
\]
be the group of $K$-rational points of $\G$. The group $G$ acts on
$Y_K$, and it is easy to see that $Y_K\to X_K$ is a $G$-Galois cover
with branch locus $(x_j)$. The ramification index of this cover at the
branch point $x_j$ is equal to $m_j$ for $j\in B\tame=B\prim\cup
B\new$ and equal to $p$ for $j\in B\wild$.  Write $Y\to
X:=\PP^1_{\Kb}$ for the pullback of $Y_K\to X_K$ to $\Kb$.

\begin{defn} \label{stabaux1def}
  The $G$-cover $f:Y\to X$ is called the {\em lift} of the special
  deformation datum $(\Zb,\omega)$ with branch locus $(x_j)$. 
\end{defn}

This definition is justified by the next proposition.

\begin{prop}  \label{stabaux1prop}
  Let $G'$ be a finite group, isomorphic to the semi-direct product
  $\ZZ/p\rtimes H$ determined by the character $\chi$.  Let $W\to X$
  be a $G'$-Galois cover. Then $W\to X$ and $Y\to X$ are isomorphic as
  mere covers if and only if the following conditions hold.
  \begin{enumerate}
  \item
    The branch locus of $W\to X$ is equal to $\{x_j\mid j\in B\}$.
  \item The cover $W\to X$ has bad reduction and the deformation datum
    occuring over the original component of the stable reduction is
    isomorphic to $(\Zb,\omega)$.
  \end{enumerate}
\end{prop}

\begin{proof} Let $W\to X$ be a $G'$-Galois cover satisfying Conditions (1)
and (2) of the proposition. After replacing $K$ by a larger field, if
necessary, we may suppose that $W\to X$ has a model $W_K\to X_K$ over
$K$ which extends to a stable model $W_R\stab\to X_R\stab$.  Let $W_R$
denote the normalization of $X_R=\PP^1_R$ in the function field of
$W_K$. We obtain a canonical map $W_R\stab\to W_R$ which is an
isomorphism except in the singular points of $\Wb$. By Condition (2),
the action of $G'$ on $\Wb$ has a kernel $I'\lhd G'$, cyclic of order
$p$. By the assumption on $G'$ and $W\to X$, we have $G'/I'\cong H$
and there exists a finite map $\Wb\to\Zb$, purely inseparable of
degree $p$ and $H$-equivariant. Clearly, the isomorphism $G'/I'\cong
H$ extends to an isomorphism $G'\cong G$.

Let $\G'$ be the schematic closure of $G'=\Aut(W_K/X_K)$ inside the
$R$-group scheme $\Aut(W_R/X_R)$. By \cite{Raynaud74}, $\G'$ is a
finite flat group scheme over $R$ with generic fiber $G'$. We claim
that the isomorphism $G'\cong G$ chosen above extends to an
isomorphism $\G'\cong\G\otimes_{W(k)}R$ of $R$-group schemes. Indeed,
it suffices to show that the closure of $I'$ inside $\G'$ is
isomorphic to $\bmu_p$ (because then $\G'$ is the minimal prolongation
of $G'$ over $R$, see \cite{Raynaud74}). Therefore, the claim follows
from Condition (2).

It is no loss of generality to identify $\G'$ with $\G\otimes_{W(k)}
R$. So $\G$ acts on $\Wb$ such that $\Zb=\Wb/\bmu_p$ and $\Xb=\Wb/\G$.
Moreover, the restriction of the map $\Wb\to\Zb$ to a sufficiently
small open subset is a $\bmu_p$-cover corresponding to the
differential $\omega$. We conclude that there exists a
$\G$-equivariant isomorphism $\Wb\iso\Yb$, see the proof of
Proposition \ref{defo1prop}. Via this isomorphism, we may regard $W_R$
as an equivariant deformation of $\Yb$. By Condition (1), the two
$\G$-covers $Y_R\to X_R$ ad $W_R\to X_R$ have the same branch locus
$(\tau_{j,R})$. Hence the proposition follows from Corollary
\ref{specialcase1cor}.  \end{proof}

\subsubsection{}

Let $(\tau_{j,R})$ be a normalized $B$-tuple over $R$,
$x_j:=\tau_{j,K}\in X:=\PP^1_{\Kb}$, and let $Y\to X$ be the lift of
$(\Zb,\omega)$ with branch locus $(x_j)$. After replacing $K$ by a
finite extension, we may assume that the natural $K$-model $Y_K\to
X_K$ extends to the stable model $Y_R\stab\to X_R\stab$. We denote by
$\Yb\stab\to\Xb\stab$ the stable reduction and by $\xb_j\in\Xb\stab$
the specialization of the branch point $x_j$.  The stable model is
obtained as a blow-up of the natural $R$-model $Y_R\to X_R$, i.e.\ we
have a commutative diagram
\begin{equation} \label{blowupdiag}
\begin{CD}
       Y_R\stab @>>> X_R\stab \\
       @VVV          @VVV     \\
       Y_R      @>>> X_R     \\
\end{CD}
\end{equation}
in which the vertical arrows are modifications which are isomorphisms
away from the critical points of $\Xb$ resp.\ the singular points of
$\Yb$. In particular, the map $\Xb\stab\to\Xb=\PP^1$ induces an
isomorphism between $\Xb$ and the original component $\Xb_0$ of
$\Xb\stab$, and contracts all other components to one of the critical
points $\tau_j$, $j\in B\prim\cup B\new$. By construction of $Y_R$,
the deformation datum over $\Xb_0\cong\Xb$ associated to
$\Yb\stab\to\Xb\stab$ is isomorphic to $(\Zb,\omega)$.

Here is the main result of this section. 

\begin{thm} \label{stabauxthm}
  The following two conditions are equivalent:
  \begin{enumerate}
  \item The curve $\Xb\stab$ consists only of the original component
    $\Xb_0$ and, for each $j\in B\tame$, a tail $\Xb_j$ which contains
    $\xb_j$ and intersects $\Xb_0\cong\Xb$ in $\tau_j$.
  \item
    For all $j\in B\new$ we have 
    \[
        \val(\tau_{j,R}-[\tau_j]) \;\geq\; 
              \frac{p\,m_j}{(p-1)\,h_j}\cdot\val(p).
    \]
    (Recall that $[\tau_j]\in W(k)\subset R$ denotes the Teichm\"uller
    lift of $\tau_j$.)
  \end{enumerate}
\end{thm}

\begin{rem}
  Under the additional assumption $B\wild=\emptyset$, the implication
  `(2) $\Rightarrow$ (1)' is already proved in \cite{special}.
\end{rem}

\begin{proof} 
The element $\tau_{j,R}-[\tau_j]\in R$ is the image
of the parameter $T_j$ (defined in \S \ref{specialcase2}) under the
classifying morphism $\Rt\to R$ of the equivariant deformation
$Y_R$. By abuse of notation, we write $T_j:=\tau_{j,R}-[\tau_j]$ as an
element of $R$. 

Let us first prove the implication `(2) $\Rightarrow$ (1)', i.e.\ we
assume that (1) holds. To find the modification $Y_R\stab\to Y_R$
which leads to the stably marked model of $Y_R$, it suffices to
consider the completion of $Y_R$ at a singular point of $\Yb$. Let
$\eta_j\in\Yb$ be one of the singular points, lying above the critical
point $\tau_j\in\Xb$ for $j\in B\tame$. Let $\xi_j\in\Zb$ be the image
of $\eta_j$. Let $\Yd_{R,j}$ denote the completion of $Y_R$ at
$\eta_j$ and $\Zd_{R,j}$ the completion of $Z_R$ at $\xi_j$. Condition
(1) of the theorem means that the germ $\Yd_{R,j}$ has potentially
good reduction.  We can write $\Zd_{R,j}=\Spec R[[z]]$. The
$\bmu_p$-torsor $\Yd_{R,j}\to \Zd_{R,j}$ is given by a Kummer equation
$y^p=u$, with $u\in R[[z]]^\times$.

Let us first treat the case $j\in B\new$. By Lemma
\ref{specialcase3lem} we may assume that
\begin{equation} \label{Kummereq2}
    y^p \;=\; u \;\equiv\; 1 \;+\; w_j\,T_j\,z^{a_j} \;+\; z^{h_j} 
         \mod{(\,p^2,\,T_j\,z^{a_j+1},\,z^{h_j+1}\,)},
\end{equation}
with $w_j\in R^\times$. We also may assume that $R$ contains an
element $\lambda$ such that $\lambda^{h_j(p-1)}=-p$. Let
$y=1+\lambda^{h_j}y'$ and $z=\lambda^pz'$. The new coordinates $y'$
and $z'$ give rise to a diagram
\[\begin{CD}
       Y_{R,j}'    @>>>  Z_{R,j}' \\
       @VVV              @VVV     \\
       \Yd_{R,j}   @>>>  \Zd_{R,j}  \\
\end{CD}\]
in which the vertical arrows are blow-ups centered in the closed point
of the germs $\Yd_{R,j}$ and $\Zd_{R,j}$, respectively. Let $\Yb_j'$
(resp.\ $\Zb_j'$) be the exceptional divisor of these blowups. The
curve $\Zb_j'$ is a projective line with parameter $z'$. To obtain the
equation for the cover $\Yb_j'\to\Zb_j'$, rewrite \eqref{Kummereq2} in
terms of $z'$ and $y'$, subtract $1$, divide by $\lambda^{ph_j}$ and
reduce modulo the maximal ideal of $R$. Using the assumption
$(p-1)h_j\val(T_j)\geq pm_j\val(p)$, one obtains an Artin-Schreier
equation
\begin{equation} \label{ArtinSchreiereq}
     {y'}^p - y' \;=\; \bar{w}_j\,{z'}^{a_j} \;+\; {z'}^{h_j} 
\end{equation}
with conductor $h_j$. It follows that $\Yb_j'$ is smooth and that
$Y_{R,j}'$ is a semistable model of the germ $\Yd_{R,j}$. In other
words, the germ $\Yd_{R,j}$ has potentially good reduction. See e.g.\ 
\cite{Raynaud98}, Lemme 4.3.3 or \cite{Lehr01}, \S 3. We conclude that
Condition (1) holds for the given index $j\in B\new$. 

Condition (1) holds for $j\in B\prim$ as well, without any assumption.
Indeed, if $j\in B\prim$ then $0<h_j<m_j$. An easy version of the
proof of Lemma \ref{specialcase3lem} shows that (for the right choice
of parameter $z$) the $\bmu_p$-torsor $Y_{R,j}\to Z_{R,j}$ is given by
a Kummer equation of the form
\[
     y^p \;=\; u \;\equiv\; 1 \;+\; \,z^{h_j} 
      \mod{(\,p^2,\,z^{h_j+1}\,)}.
\]
After a change of coordinate $y=1+\lambda^{h_j}y'$ and
$z'=\lambda^pz'$ we again find an Artin--Schreier equation with
conductor $h_j$, which shows that the germ $\Yd_{R,j}$ has potentially
good reduction. This finishes the proof of the implication `(2)
$\Rightarrow$ (1)'.

It remains to prove the implication `(1) $\Rightarrow$ (2)'. Suppose
that Condition (2) is violated for some $j\in B\new$. We may also
assume that $R$ contains an element $\lambda$ with
$\lambda^{pm_j}=T_j$. After substituting $y=1+\lambda^{h_j}y'$ and
$z=\lambda^pz'$ into \eqref{Kummereq2} and reducing modulo the maximal
ideal of $R$, we find that the cover $\Yb_j'\to\Zb_j'$ has equation
\[
    {y'}^p \;=\; \bar{w}_j\,{z'}^{a_j} \;+\; {z'}^{h_j}.
\]
It follows that the cover $\Yb_j'\to\Zb_j'$ is an $\balpha_p$-torsor
(away from the point $z'=\infty$) corresponding to the exact
differential $\omega'=\diff(\bar{w}_j{z'}^{a_j}+{z'}^{h_j})$.
The differential $\omega'$ has a zero of order $a_j-1$ at $z'=0$ and
$m_j$ further simple zeros. Above each point
where $\omega'$ has a zero, $\Yb_j'$ has a singularity. We conclude
that the germ $\Yd_{R,j}$ does not have potentially good reduction. In
other words, Condition (1) of Theorem \ref{stabauxthm} does not hold.
This finishes the proof of Theorem \ref{stabauxthm}.  \end{proof}

The proof shows that, if both the conditions of the theorem hold, then
the closed rigid disk
\[
     D_j \;:=\; \{\;x \,\mid\, \val(x-[\tau_j])\geq 
         \frac{p\,m_j}{(p-1)\,h_j}\cdot\val(p)\;\} 
               \;\subset\; X_K^{\rm rig}
\]
is equal to the set of points on $X_K$ which specialize to the tail
$\Xb_j\subset X_R\stab$. We remark at this point that, since
$\frac{pm_j}{(p-1)h_j}<1$ for $j\in B\new$, the family of disks
$(D_j)_{j\in B}$ is independent of the choice of $(\tau_{j,R})$, as
long as Condition (2) of Theorem \ref{stabauxthm} holds.

\subsubsection{} \label{auxmon}

Let us now assume that the branch points $x_j$ of the lift $f:Y\to X$
of $(\Zb,\omega)$ are all $K_0$-rational. This is equivalent to the
condition that $\tau_{j,R}\in W(k)$. Then Condition (2) of Theorem
\ref{stabauxthm} is automatically verified. Moreover, the cover
$f:Y\to X$ has a natural model $f_{K_0}:Y_{K_0}\to
X_{K_0}=\PP^1_{K_0}$ over $K_0$ (which is {\em not} Galois, since
$K_0$ does not contain the $p$th roots of unity).  Let
\[
     \kappab\ab:\Gal(\Kb/K_0) \To \Aut_k(\Yb\stab)
\]
be the absolute monodromy action associated to the $K_0$-model
$f_{K_0}$ of $f$, see \S \ref{monodromy}. The goal of this subsection
is to describe this action in detail.

First some notation. Let $\Yb_0$ be the unique irreducible component
of $\Yb\stab$ which lies above the original component
$\Xb_0\subset\Xb\stab$. We may identify $\Yb_0$ with the normalization
of $\Yb$ or, what amounts to the same thing, the smooth projective
model of the function field $k(\Zb)^{1/p}$. Let us choose, for all
$j\in B\tame$, a point $\eta_j\in\Yb_0$ above $\tau_j\in\Xb$. We
denote by $\xi_j\in\Zb$ the image of $\eta_j$, by $H_j\subset H$ the
stabilizer of $\xi_j$ and by $I\subset G$ the cyclic subgroup of order
$p$. Set $G_j:=I\cdot H_j\subset G$. Let $\Yb_j\subset\Yb\stab$
(resp.\ $\Xb_j\subset\Xb\stab$) be the component which intersects
$\Yb_0$ in $\eta_j$ (resp.\ intersects $\Xb$ in $\tau_j$). By Theorem
\ref{stabauxthm}, $\Xb\stab=\Xb\cup\cup_j\Xb_j$. Furthermore, the
proof of the theorem shows that $(\Yb_j\to\Xb_j,\eta_j,H_j)$ is a
pointed $G_j$-tail cover with ramification type $(m_j,h_j)$ which is
totally ramified above $\tau_j\in\Xb_j$, see \S \ref{Gdefdat2}.

Let us denote by $\Aut_k^0(\Yb\stab)$ the group of $k$-linear
automorphisms of $\Yb\stab$ which normalize the action of $G$, commute
with the action of $H\subset G$ and induce the identity on $\Yb_0$.
Similarly, let $\Aut_k^0(\Yb_j)$ denote the group of (outer)
automorphisms of the pointed $G_j$-tail cover $\fb_j$. Recall that
these are $k$-linear automorphisms of $\Yb_j$ which normalize the
action of $G_j$, commute with the action of $H_j\subset G_j$ and fix
the point $\eta_j$. For each $j$ we have a natural homomorphism
\begin{equation} \label{auxmoneq1}
     \Aut_k^0(\Yb\stab) \;\To\; \Aut_k^0(\Yb_j).
\end{equation}
Let $\alpha$ be a generator of $I\lhd G$. For
$\vphi\in\Aut_k^0(\Yb)$ we have
$\vphi\circ\alpha\circ\vphi^{-1}=\alpha^{\rho}$ for some $\rho\in
\FF_p^\times$; this defines a character
\begin{equation} \label{auxmoneq2}
    \rho:\Aut_k^0(\Yb) \;\To\; \FF_p^\times.
\end{equation}
Similarly, one defines characters
$\rho_j:\Aut_k^0(\Yb_j)\to\FF_p^\times$ for all $j\in B\tame$. The
composition of $\rho$ with $\kappab\ab$ is the cyclotomic character
modulo $p$. The kernel of the character $\rho$ is the subgroup
$\Aut_G^0(\Yb\stab)\subset\Aut_k^0(\Yb\stab)$ of automorphisms which
commute with the $G$-action. Compare with \S \ref{Gdefdat4}.

\begin{prop} \label{auxmonprop}
\begin{enumerate}
\item
  We have a canonical isomorphism
  \[
    \Aut_k^0(\Yb\stab) \;\cong\;
    \{\,(\rho;\vphi_j)\in\FF_p^\times\times\prod_j \Aut_k^0(\Yb_j) 
      \;\mid\; \rho=\rho_j(\vphi_j)\; \forall j\;\}.
  \] 
\item
  The group $\Aut_k^0(\Yb_j)$ is cyclic of order
  $(p-1)h_j$. Therefore, the group $\Aut_k^0(\Yb\stab)$ has order
  $(p-1)\prod_jh_j$.   
\item
  The restriction of $\kappab\ab$ to the component $\Yb_j$ is a
  surjective homomorphism
  \[
     \kappab\ab_j: \Gal(\Kb/K_0) \;\To\; \Aut_k^0(\Yb_j).
  \]
\item
  The image of $\kappab\ab$ is a cyclic subgroup of $\Aut_k^0(\Yb\stab)$,
  of order
  \[
          N \;:=\; (p-1)\cdot {\rm lcm}_j (h_j).
  \]
\end{enumerate}
\end{prop}
   
\begin{proof} Part (1) is left to the reader. We also leave it to the reader
to check that the image of the monodromy action $\kappab\ab$ is
contained in $\Aut_k^0(\Yb\stab)$.  

The explicit description of the change of coordinates
$y=1+\lambda^{h_j}y'$, $z=\lambda^pz'$ in the proof of Theorem
\ref{stabauxthm} shows that the restriction of $\kappab\ab$ to $\Yb_j$
is cyclic of order $(p-1)h_j$ (it corresponds to the cyclic Galois
extension $K_0(\lambda)/K_0$ where $\lambda^{(p-1)h_j}=-p$). In
particular, the group $\Aut_k^0(\Yb_j)$ has a cyclic subgroup of order
$(p-1)h_j$.  Therefore, in order to prove the proposition it suffices
to show that the group $\Aut_k^0(\Yb_j)$ has exact order $(p-1)h_j$.

First, we claim that $\Aut_k^0(\Yb_j)$ is cyclic. Indeed, the natural
action of $\Aut_k^0(\Yb_j)$ on $\Yb_j$ descents to a faithful action
on the quotient $\Xb_j=\Yb_j/G_j$ which fixes the two points $\tau_j$
and $\xb_j$. But $\Xb_j$ has genus $0$, so the claim follows. Now
assume that $\Aut_k^0(\Yb_j)$ has order $(p-1)h_jn$. Let
$\vphi\in\Aut_k^0(\Yb_j)$ be an element of order $h_jn$. Since $\vphi$
lies in the kernel of the character $\rho_j$ it commutes with the
action of $G_j$ on $\Yb_j$. Let
$\Yb_j':=\Yb_j/\gen{\vphi}\to\Xb_j':=\Xb_j/\gen{\vphi}$ be the induced
$G_j$-cover. Recall that the action of $I\cong\ZZ/p$ on $\Yb_j$ has
conductor $h_j$ at the point $\eta_j$. Let $\eta_j'$ denote the image
of $\eta_j$ on $\Yb_j'$. A classical calculation (see e.g.\
\cite{Raynaud98}, Lemme 1.1.2) shows that the action of $I$ on
$\Yb_j'$ has conductor $1/n$ at $\eta_j'$. We conclude that
$n=1$. This finishes the proof of the proposition.  \end{proof}

The proposition implies that the $G$-cover $f:Y\to X$ has stable
reduction over the unique tame extension $K\stab/K_0$ of degree $N$.


\section{Construction of three point covers by lifting}  \label{main}

Let $\fb:\Yb\to\Xb$ be a special $G$-map, see Definition
\ref{specialGmapdef}. A {\em lift} of $\fb$ is a three point
$G$-cover with stable reduction $\fb$. Let $L(\fb)$ denote the set of
isomorphism classes of all lifts of $\fb$. The Galois group
$\Gal(\Kb/K_0)$ acts on $L(\fb)$ in a natural way. In this section we
determine the structure of $L(\fb)$ as a $\Gal(\Kb/K_0)$-set. In
particular, we prove that $L(\fb)$ is nonempty and that the action of
$\Gal(\Kb/K_0)$ is tamely ramified.

At the end, we derive various consequences of the above result, concerning the
field of moduli of three point covers, and we discuss a aprticular class of
examples, called {\em genus zero dessins of prime degree}.

\subsection{The set of lifts of $\fb$} \label{lifts}

\subsubsection{} \label{lifts1}

Let $k$ be an algebraic closure of the finite field $\FF_p$,
$W(k)$ the ring of Witt vectors over $k$ and $K_0$ the fraction
field of $W(k)$. We fix an algebraic closure $\Kb$ of $K_0$. In this
section, $K/K_0$ always denotes a finite field extension contained in
$\Kb$ and $R$ denotes the ring of integers of $K$. Let
$(\Zb_0,\omega_0;\Yb_j)$ be a special $G$-deformation datum and let
$\fb:\Yb\to\Xb$ be the associated special $G$-map, see \S
\ref{Gdefdat}.

\begin{defn}
  A {\em lift} of $\fb$ over $R$ is a stable $G$-map $f_R:Y_R\to
  X_R$ over $R$ (see Definition \ref{stabGmapdef}), together with an
  isomorphism $\fb\cong f_R\otimes_R k$ of stable $G$-maps and an
  isomorphism $X_R\otimes_R K\iso(\PP^1_K,\{0,1,\infty\})$ of marked
  curves. Given a lift $f_R:Y_R\to X_R$, we will usually identify
  $\fb$ with the special fiber of $f_R$ and $\PP^1_K$ with the general
  fiber of $X_R$. Two lifts $f_R:Y_R\to X_R$ and $f_R':Y_R'\to
  X_R'$ over $R$ are said to be isomorphic if there exists an
  isomorphism $f_R\iso f_R'$ of stable $G$-maps which induces the
  identity on the special fiber and on $X_K=X_K'=\PP^1_K$.
\end{defn}

If $f_R:Y_R\to X_R$ is a lift of $\fb$, then the general fiber
$f_K:Y_K\to X_K=\PP^1_K$ of $f_R$ is a three point $G$-cover and $f_R$
is its stable model. In particular, $f_K$ has mild reduction and $\fb$
is the stable reduction of $f_K$. If $f_R$ and $f_R'$
are isomorphic lifts of $\fb$ then the isomorphism $f_R\iso f_R'$
is unique.  (This is a general property of deformations of stably
marked curves, see \cite{Knudsen83}.)

\subsubsection{} \label{lifts2}

Let $L(\fb,K)$ denote the set of isomorphism classes of lifts of
$\fb$ over $R$, the ring of integers of $K$. If $K'/K$ is a finite
extension then we have a natural inclusion $L(\fb,K)\inj L(\fb,K')$.
Let 
\[
         L(\fb) \;:=\; \cup_K L(\fb,K)
\]
be the union, taken over all finite extensions $K/K_0$
contained in $\Kb$. It is clear that $L(\fb)$ is a finite set. A
priori, $L(\fb)$ may be empty; we will see later that it is not.

There is a natural continuous action of $\Gal(\Kb/K_0)$ on $L(\fb)$
from the left, defined as follows. Let $f_R:Y_R\to X_R$ be a lift
of $\fb$ over $R$, the ring of integers of $K$, and let
$\sigma\in\Gal(\Kb/K_0)$. We may assume that $K/K_0$ is Galois. Let
$\op{\sigma}{f_R}:\op{\sigma}{Y_R}\to\op{\sigma}{X_R}$ be the pullback
of $f_R$ by the automorphism $\tilde{\sigma}:\Spec R\iso\Spec R$
induced by $\sigma$. Since $\sigma$ induces the identity on the
residue field $k$, we have a canonical isomorphism $f_R\otimes_Rk\iso
\op{\sigma}{f_R}\otimes_R k$. We may therefore consider
$\op{\sigma}{f_R}$ as a deformation of $\fb$. This defines the action
of $\Gal(\Kb/K_0)$ on $L(\fb)$. 

Recall that $\Aut_G^0(\fb)$ was defined in \S \ref{Gdefdat4} as the
group of automorphisms of $\fb$ which induce the identity on the
original component. There is a natural right action of $\Aut_G^0(\fb)$
on $L(\fb)$, defined as follows. If $\gamma\in\Aut_G^0(\fb)$ and
$f_R:Y_R\to X_R$ is a lift of $\fb$, then we define $f_R^\gamma$ as
the stable $G$-map $f_R:Y_R\to X_R$ together with the isomorphism
$\Yb\iso Y_R\otimes_Rk$ which is the composition with
$\gamma:\Yb\iso\Yb$ with the canonical isomorphism $\Yb\iso
Y_R\otimes_R k$. This defines an action of $\Aut_G^0(\fb)$ on $L(\fb)$
from the right which commutes with the action of $\Gal(\Kb/K_0)$.  We
define
\[
     \tilde{L}(\fb) \;:=\; L(\fb)/\Aut_G^0(\fb).
\]
It is clear that elements of $\tilde{L}(\fb)$ correspond naturally to
isomorphism classes of three point $G$-covers $f:Y\to X$ whose
stable reduction is isomorphic to $\fb$. By definition we have:

\begin{prop}
  Let $f_R:Y_R\to X_R$ be a lift of $\fb$ and denote by $f:Y\to X$
  the geometric generic fiber of $f_R$. 
  \begin{enumerate}
  \item
    The stabilizer in $\Gal(\Kb/K_0)$ of $f_R$, considered as an
    element of $L(\fb)$, is equal to $\Gamma_f\stab$, as defined in \S
    \ref{monodromy}. 
  \item 
    The stabilizer in $\Gal(\Kb/K_0)$ of $f_R$, considered as an
    element of $\tilde{L}(\fb)$, is equal to $\Gamma_f\inn$, as
    defined in \S \ref{monodromy}.
  \end{enumerate}
\end{prop}

\subsubsection{The auxiliary cover} \label{patch1}

Let $f_R:Y_R\to X_R$ be a lift of $\fb$. For $j\in B\prim\cup B\new$,
let $D_j\subset X_K^{\rm rig}$ denote the closed rigid disk
corresponding to the tail $\Xb_j$. By definition, a $K$-rational point
$x\in D_j$ specializes to a point $\xb\in\Xb_j$. In particular, for
$j\in B\prim$ the branch point $x_j\in\{0,1,\infty\}$ of $f$ lies on
$D_j$. Note that the model $X_R$ of $X_K=\PP^1_K$ is uniquely
determined by the family of disks $(D_j)_j$. A priori, the disks $D_j$
(and therefore the model $X_R$) may depend on the lift $f_R$. We are
now going to show that they are in fact independent of $f_R$. This is
an important steps towards our main result.

Let us choose, for all $j\in B\new$, a $K$-rational point $x_j\in
D_j$. We call these points the {\em auxiliary branch points}.  By
\cite{Raynaud98}, \S 3.2, the choice of the auxiliary branch points
yields a certain $G_0$-cover $f\aux:Y\aux\to X$, called the {\em
auxiliary cover} of $f$. By definition, $f\aux$ has branch locus
$\{x_j\mid j\in B\}$ and the same type of bad reduction as $f$. To be
more precise, let $\hat{X}_R$ denote the formal completion of $X_R$
along $\Xb$ and let $\D_j\subset\hat{X}_R$, for $j\in B\prim\cup
B\new$, denote the open subset corresponding to
$\Xb_j-\{\tau_j\}\subset\Xb$, i.e.\ $\D_j$ is the closed formal disk
whose generic fiber is $D_j$. Let $\U^\circ\subset\hat{X}_R$ denote
the open subset corresponding to $\Xb-\cup_j\Xb_j\subset\Xb$.
Finally, let $\X_j$ denote the completion of $X_R$ at $\tau_j\in\Xb$
and set $\U:=\U^\circ\cup\cup_j\X_j$. Recall that part of our special
$G$-deformation datum $(\Zb_0,\omega_0;\Yb_j)$ is a subgroup
$G_0=I_0\cdot H_0$. The auxiliary cover of the lift $f_R:Y_R\to X_R$
is a $G_0$-cover $f\aux:Y\aux\to X$, characterized by the following
properties, see \cite{Raynaud98}, \S 3.2.
\begin{itemize}
\item[(a)]
  The cover $f\aux$ has a stable model $f_R\aux:Y_R\aux\to X_R\aux$
  over $R$. Moreover, we have a (unique) isomorphism $X_R\cong X_R\aux$
  extending the identity on the generic fiber.
\item[(b)] For $j\not\in B\wild$, the restriction of $f_R\aux$ to the
  formal disk $\D_j$ is a tame cover of formal schemes, ramified along
  the closure of $x_j$ in $\D_j$.
\item[(c)]
  There exists a $G$-equivariant isomorphism
  \begin{equation} \label{patcheq1}
     \hat{Y}_R \times_{\hat{X}_R}\U \;\cong\; 
          \Ind_{G_0}^G(\hat{Y}_R\aux\times_{\hat{X}_R}\U)
  \end{equation}
  of formal $\U$-schemes. 
\end{itemize}
It follows that the branch locus of $f\aux:Y\aux\to X$ is precisely
the set of points $\{x_j\mid j\in B\}$.

\begin{prop} \label{patchprop1}
  \begin{enumerate}
  \item
    The cover $f\aux$ is the (unique) lift of $(\Zb_0,\omega_0)$
    with branch points $(x_j)$, see \S \ref{stabaux1}. In
    particular, as a mere cover, $f\aux$ does not depend on the lift
    $f_R$ of $\fb$.
  \item 
    The disks $D_j$ are independent of the lift $f_R:Y_R\to
    X_R$. Moreover, each disk $D_j$ contains a $K_0$-rational point.
  \end{enumerate}
\end{prop}

\begin{proof} It follows immediately from Property (c) that the deformation
datum over the original component associated to the stable reduction
of $f\aux$ can be identified with $(\Zb_0,\omega_0)$. Hence (1)
follows from Proposition \ref{stabaux1prop}.  By Property (a), the
cover $f\aux$ satisfies Condition (1) of Theorem
\ref{stabauxthm}. Therefore, (2) follows from Theorem
\ref{stabauxthm}.  \end{proof}

We remark that the proof of Proposition \ref{patchprop1} is the only
place where we use the subtle implication `(1)$\Rightarrow$(2)' of
Theorem \ref{stabauxthm}.

\subsection{Patching data} \label{patch}

Lifting the special $G$-map $\fb$ to a three point $G$-cover $f:Y\to
X$ can be done in three steps. The first step consists in lifting the
tail covers $\fb_j:\Yb_j\to\Xb_j$ to Galois covers $\E_j\to\D_j$ of
the formal disks $\D_j$ corresponding to $D_j\subset X$. Since the
restriction of $\fb_j$ to $\Xb_j-\{\tau_j\}$ is tame, such a lift
exists and is unique.  The second step consists in lifting the
(inseparable) cover $\Yb_0\to\Xb_0$ to a Galois cover $\V\to\U$, where
$\U$ is the `formal neighborhood' of
$\hat{X}_R-\cup_j\D_j$ defined in \S \ref{patch1}. The existence
and uniqueness of such a lift follows from the results of \S
\ref{defo}. In the third step we patch the covers $\E_j\to\D_j$ and
$\V\to\U$ together along the boundaries $\B_j$ of the disks $\D_j$ and
obtain the desired $G$-cover $f:Y\to X$. Unlike the first two steps,
the third step is not unique; it depends on the choice of what we
call a {\em patching datum} for $\fb$. The construction we have just
sketched amounts to the existence of a bijection
\[
        L(\fb) \;\liso\; P(\fb)
\]
between the set of lifts of $\fb$ and the set $P(\fb)$ of all
patching data. Essentially by construction, this bijection is
$\Gal(\Kb/K_0)$-invariant. Moreover, the results of \S \ref{auxmon}
can be used to determine the complete structure of $P(\fb)$ as
$\Gal(\Kb/K_0)$-set. Intuitively, the Galois group $\Gal(\Kb/K_0)$
acts on $P(\fb)$ by `Dehn-twists' along the boundaries $\B_j$.

\subsubsection{} \label{patch2}

For each index $j\in B$, let us now choose, once and for all, a
$K_0$-rational point $x_j\in X$ such that the closure of $x_j$
in $X_0=\PP^1_{W(k)}$ intersects the special fiber $\Xb_0=\PP^1_k$ in
the critical point $\tau_j$. As usual, we assume that $x_j\in
\{0,1,\infty\}$ for $j\in B_0$. Let $f\dagg:Y\dagg\to X=\PP^1_{\Kb}$
be the lift of the special deformation datum $(\Zb_0,\omega_0)$
with branch points $(x_j)_{j\in B}$, see \S \ref{stabaux}. By
definition, $f\dagg$ is a $G\dagg$-cover, where
$G\dagg:=\bmu_p(\Kb)\rtimes H_0$.  If $f_R$ is a lift of $\fb$,
then Proposition \ref{patchprop1} shows that $f\dagg$ is isomorphic,
as a mere cover, to the auxiliary cover $f\aux$ associated to $f_R$.
In the following, we will not identify $f\dagg$ with $f\aux$, because
we want to emphasize the fact that such an identification is not
canonical.

Let $K/K_0$ be a finite extension over which $f\dagg$ has stable
reduction and let $f_R\dagg:Y_R\dagg\to X_R\dagg$ (resp.\
$\fb\dagg:\Yb\dagg\to\Xb\dagg$) denote the stable model (resp.\ the
stable reduction) of $f\dagg$. By Theorem \ref{stabauxthm}, the curve
$\Xb\dagg$ consists of the original component $\Xb_0=\PP^1_k$ and the
tails $\Xb_j\dagg$, where $j$ runs over the set $B\prim\cup B\new$.
Moreover, the deformation datum over the original component associated
to $\fb\dagg$ is canonically isomorphic to $(\Zb_0,\omega_0)$.  We may
and will identify the (unique) component $\Yb_0\dagg$ of $\Yb\dagg$
over the original component with the distinguished component $\Yb_0$
of $\Yb$ (by definition, $\Yb_0=\Zb_0^{1/p}$, see \S \ref{Gdefdat3}).
Recall that we have distinguished points $\eta_j\in\Yb_0$ above
$\tau_j\in\Xb_0$; they are part of the $G$-deformation datum
$(\Zb_0,\omega_0;\Yb_j)$. For $j\not\in B\wild$ we let $\Yb_j\dagg$
denote the irreducible component of $\Yb\dagg$ which intersects
$\Yb_0$ in $\eta_j$. Write $\eta_j\dagg\in\Yb_j\dagg$ for the point
where $\Yb_j\dagg$ intersects $\Yb_0$. Let $G_j\dagg\subset G\dagg$ be
the subgroup generated by $\bmu_p(\Kb)$ and $H_j$. One checks that
$(\fb_j\dagg:\Yb_j\dagg\to\Xb_j\dagg,\eta_j\dagg,H_j)$ is a pointed
$G_j\dagg$-tail cover with ramification type $(m_j,h_j)$, totally
ramified above $\tau_j$. See \S \ref{Gdefdat2}.

Let $\Aut_k^0(\Yb\dagg)$ denote the group of $k$-linear automorphisms
of $\Yb\dagg$ which normalize the action of $G\dagg$, commute with the
action of $H_0\subset G\dagg$ and restrict to the identity on $\Yb_0$.
We have seen in \S \ref{auxmon} that the absolute monodromy action on
$\fb\dagg$ is a homomorphism
\[
     \kappab\dagg:\Gal(\Kb/K_0) \;\To\; \Aut_k^0(\Yb\dagg).
\]
Similarly, let $\Aut_k^0(\Yb_j\dagg)$ denote the group of $k$-linear
automorphisms of $\Yb_j\dagg$ which normalize the action of
$G_j\dagg$, commute with the action of $H_j\subset G_j\dagg$ and fix
the point $\eta_j$. By Proposition \ref{auxmonprop}, the group
$\Aut_k^0(\Yb_j\dagg)$ is cyclic of order $(p-1)h_j$; furthermore,
restricting $\kappab\dagg$ to $\Yb_j\dagg$ yields a surjective
homomorphism
\[
     \kappab_j\dagg:\Gal(\Kb/K_0) \;\To\;
       \Aut_k^0(\Yb_j\dagg) \;\cong\; \ZZ/(p-1)h_j.  
\]

\subsubsection{}

Fix $j\not\in B\wild$. We denote by $\Xb_{j,\infb}$ (resp.\
$\Xb_{j,\infb}\dagg$) the generic point of the completion of the tail
$\Xb_j$ (resp.\ the tail $\Xb_j\dagg$) at the point $\tau_j$. We set
$\Yb_{j,\infb}:=\Yb_j\times_{\Xb_j}\Xb_{j,\infb}$ and
$\Yb_{j,\infb}\dagg:=\Yb_j\dagg\times_{\Xb_j\dagg}\Xb_{j,\infb}\dagg$.
Let $P_j(\fb)$ be the set of pairs $(\lambda,\vphi_{\infb})$
consisting of the following objects:
\begin{itemize}
\item 
  $\lambda:G\dagg\iso G_0$ is an isomorphism which induces the
  identity on $H_0$, and  
\item $\vphi_{\infb}:\Ind_{G_j}^G(\Yb_{j,\infb})\iso
  \Ind_{G_j\dagg,\lambda}^G(\Yb_{j,\infb}\dagg)$ is a $G$-equivariant
  isomorphism such that the following holds.
  \begin{enumerate}
  \item
    We have $\vphi_{\infb}(\eta_j)=\eta_j\dagg$.
  \item
    The induced isomorphism
    $\psi_{\infb}:\Xb_{j,\infb}\iso\Xb_{j,\infb}\dagg$ extends to
    an isomorphism $\psi:\Xb_j\iso\Xb_j\dagg$.
  \item
    If $j\in B\prim$ then $\psi(\xb_j)=\xb_j\dagg$.
  \end{enumerate}
  (The notation `$\Ind_{G_j\dagg,\lambda}^G$' means that we embed
  $G_j\dagg$ into $G$ via $\lambda:G\dagg\iso G_0\subset G$.)
\end{itemize}
The group $\Aut_k^0(\Yb_j\dagg)$ acts on the set $P_j(\fb)$ from the
left, as follows. Let $(\lambda,\vphi_{\infb})$ be an element of
$P_j(\fb)$ and $\gamma$ an element of $\Aut_k^0(\Yb_j\dagg)$.  Let
${\rm inn}(\gamma^{-1}):G\dagg\iso G\dagg$ denote the unique
automorphism of $G\dagg$ which sends an element $\alpha\in I\dagg$ to
$\gamma^{-1}\circ\alpha\circ\gamma$ and induces the identity on
$H_0\subset G\dagg$. Set $\lambda':=\lambda\circ{\rm
  inn}(\gamma^{-1})$. One checks that $\gamma$ extends uniquely to a
$G$-equivariant isomorphism
\[
    \tilde{\gamma}: \Ind_{G_j\dagg,\lambda}^G(\Yb_j\dagg)
        \;\liso\; \Ind_{G_j\dagg,\lambda'}^G(\Yb_j\dagg).
\]
Let $\vphi_{\infb}'$ be the composition of $\vphi_{\infb}$ with the
local isomorphism induced by $\tilde{\gamma}$. Then
$(\lambda',\vphi_{\infb}')$ is an element of $P_j(\fb)$. This defines
the action of $\Aut_j^0(\Yb_j\dagg)$ on $P_j(\fb)$. It is easy to see
that this action makes $P_j(\fb)$ a principal homogeneous
$\Aut_k^0(\Yb_j\dagg)$-space. Via the monodromy action
$\kappab_j\dagg$, we obtain a transitive Galois action of
$\Gal(\Kb/K_0)$ on $P_j(\fb)$.

Define
\[
      P(\fb) \;:=\; \{\,(\lambda;\vphi_{j,\infb}) \;\mid\;
           (\lambda,\vphi_{j,\infb})\in P_j(\fb)\;\forall j\;\}.
\]
An element of $P(\fb)$ is called a
{\em patching datum}. The actions of the groups $\Aut_k^0(\Yb_j\dagg)$ on
$P_j(\fb)$ extend uniquely to an action of $\Aut_k^0(\Yb\dagg)$ on
$P(\fb)$ (use Proposition \ref{auxmonprop} (1)). Again, this action
is transitive and fixed point free. We consider $P(\fb)$ as a
$\Gal(\Kb/K_0)$-set, via the monodromy action $\kappab\dagg$.  

\begin{prop} \label{patchprop2}
   The set $P(\fb)$ of patching data has exactly $(p-1)\cdot\prod_jh_j$
   elements. The $\Gal(\Kb/K_0)$-orbits have length
   \[
         N \;:=\; (p-1)\cdot{\rm lcm}_j(h_j).
   \]
   In particular, the action of $\Gal(\Kb/K_0)$ on $P(\fb)$ is tamely
   ramified. 
\end{prop}

\begin{proof}
By Proposition \ref{auxmonprop} and the freeness of the action of 
$\Aut_k^0(\Yb\dagg)$ on $P(\fb)$, it suffices to show that $P(\fb)$ is
nonempty. But this follows from Lemma \ref{Gtaillem}. 
\end{proof}

\subsubsection{}

We shall now define a $\Gal(\Kb/K_0)$-equivariant map
\begin{equation} \label{patcheq2}
         L(\fb) \To P(\fb),
\end{equation}
as follows. Let $f_R:Y_R\to X_R$ be a lift of $\fb$. Let
$f\aux:Y\aux\to X$ be the auxiliary cover associated to $f$ and the
set $(x_j)_{j\in B}$. Let $f_R\aux:Y_R\aux\to X_R$ denote the stable
model of the auxiliary cover $f\aux$. By Proposition \ref{patchprop1}
(1), there exists a $K$-linear isomorphism $\vphi_K:Y_K\aux\iso
Y_K\dagg$ which induces the identity on $X_K$. By the uniqueness of
the stable model, $\vphi_K$ extends to an isomorphism
$\vphi_R:Y_R\aux\iso Y_R\dagg$. Let $\vphib:\Yb\aux\iso\Yb\dagg$
denote the restriction of $\vphi_R$ to the special fiber. After
composing $\vphi_K$ with an element of $G\dagg$, we may assume that
$\vphi_K$ commutes with the action of $H_0$ (which is a subgroup of
$G_0$ and $G\dagg$) and that $\vphib$ induces the identity on $\Yb_0$
(which is a component of $\Yb\aux$ and $\Yb\dagg$). Note that this
determines $\vphi_K$ uniquely. Let $\lambda:G_0\iso G\dagg$ be the
isomorphism $\alpha\mapsto \vphi_K\circ\alpha\circ\vphi_K^{-1}$. Fix
$j\not\in B\wild$. Let $\Yb_j\aux$ denote the component of
$\Yb\aux$ which intersects $\Yb_0$ in $\eta_j$, and set
$\Yb_{j,\infb}\aux:=\Yb_j\aux\times_{\Xb_j}\Xb_{j,\infb}$. The
isomorphism \eqref{patcheq1} induces a $G$-equivariant isomorphism
\begin{equation} \label{patcheq3}
   \Ind_{G_j}^G(\Yb_{j,\infb}) \;\liso\;
              \Ind_{G_0}^G(\Yb_{j,\infb}\aux).
\end{equation}
Let $\vphi_{j,\infb}:\Ind_{G_j}^G(\Yb_{j,\infb})\iso
\Ind_{G\dagg}^G(\Yb_{j,\infb}\dagg)$ be the composition of this
isomorphism with the isomorphism induced by $\vphib$. By construction,
the tuple $(\lambda;\vphi_{j,\infb})$ is a patching datum. This
defines the map \eqref{patcheq2}. 

We claim that \eqref{patcheq2} is $\Gal(\Kb/K_0)$-equivariant. Indeed,
for any $\sigma\in\Gal(\Kb/K_0)$, we may identify the auxiliary cover
of the twisted cover $\op{\sigma}{f_R}:=f_R\otimes_R^\sigma R$ with
the $\sigma$-twist $\op{\sigma}{(f_R\aux)}:=f_R\aux\otimes_R^\sigma R$
of the auxiliary cover of $f_R$. Let
$\kappa_{R,\sigma}\dagg:\op{\sigma}{(Y_R\dagg)}\iso Y_R\dagg$ be the
unique extension of the canonical $K$-linear isomorphism
$\op{\sigma}{Y_K\dagg}\iso Y_K\dagg$ induced from the $K_0$-model
$Y_{K_0}\dagg$. The restriction of $\kappa_{R,\sigma}$ to the special
fiber is equal to $\kappab\dagg(\sigma)$, i.e.\ the monodromy action
of $\sigma$ on $\Yb\dagg$. The composition
$\vphi_R':=\kappa_{R,\sigma}\dagg\circ\vphi_R$ is the unique
isomorphism which induces an $H_0$-equivariant isomorphism of mere
covers $f_K\aux\cong f_K\dagg$ and whose restriction to the special
fiber is the identity on $\Yb_0$. It follows that the map
\eqref{patcheq2} maps the class of the twisted lift
$\op{\sigma}{f_R}$ to the image of $(\lambda;\vphi_{j,\infb})$ under
the action of $\kappab\dagg(\sigma)\in\Aut_k^0(\Yb\dagg)$. This proves
the claim that \eqref{patcheq2} is $\Gal(\Kb/K_0)$-equivariant.

\begin{thm} \label{liftthm}
  The map \eqref{patcheq2} is an isomorphism of $\Gal(\Kb/K_0)$-sets. 
\end{thm}

\begin{proof} We will show that \eqref{patcheq2} is an isomorphism by
constructing an inverse $P(\fb)\to L(\fb)$. Let
$(\lambda;\vphib_{j,\infb})\in P(\fb)$ be a patching datum. We have to
construct a lift $f_R:Y_R\to X_R$ of $\fb$ which is mapped to
$(\lambda;\vphib_{j,\infb})$ by \eqref{patcheq2}.

By definition, the local isomorphism
$\psib_{j,\infb}:\Xb_{j,\infb}\iso\Xb_{j,\infb}\dagg$ extends to a
global isomorphism $\psib_j:\Xb_j\iso\Xb_j\dagg$. Moreover, if $j\in
B\prim$ then $\psi_j(\xb_j)=\xb_j\dagg$. Let $\psib:\Xb\iso\Xb\dagg$
denote the isomorphism which is equal to the identity on the original
component and which is equal to $\psib_j$ on the tails.  We define
$X_R:=X_R\dagg$ as a marked curve over $R$. However, we identify the
special fiber of $X_R$ with $\Xb$ via the isomorphism $\psib$.

Let $\hat{X}_R$ denote the formal completion of $X_R$ along $\Xb$. For
$j\not\in B\wild$, let $\D_j\subset\hat{X}_R$ be the open formal
subscheme corresponding to $\Xb_j-\{\tau_j\}\subset\Xb$ and let $\X_j$
be the completion of $X_R$ at $\tau_j$. Let $\U^\circ\subset\hat{X}_R$
be the open formal subscheme corresponding to $\Xb_0-\{\tau_j\mid j\in
B\prim\cup B\new\}\subset\Xb$. Let $\U$ denote the disjoint union of
$\U^\circ$ with all the $\X_j$. The fiber product
\[
    \B_j \;:=\; \U\times_{\hat{X}_R}\D_j 
       \;\cong\; \Spec R[[T]]\{T^{-1}\}
\]
is the `boundary' of the formal disk $\D_j$. Its special fiber
$\B_j\otimes_R k$ is canonically isomorphic to $\Xb_{j,\infb}\dagg$
and isomorphic to $\Xb_{j,\infb}$ via $\psib_{j,\infb}$.  Define
\[
      \V \;:=\; \Ind_{G\dagg,\lambda}^G
           (\hat{Y}_R\dagg\times_{\hat{X}_R}\U).
\]
Furthermore, let $\E_j\to\D_j$ be the (unique) tame lift of the
(not necessarily connected) $G$-cover 
\[
     \Yb\times_{\Xb}(\Xb_j-\{\tau_j\}) \;\To\; \Xb_j-\{\tau_j\}
         \;\lpfeil{\psib_j}\; \Xb_j\dagg-\{\tau_j\}
\]
which is unramified away from the section $x_{j,R}:\Spec
R\to\D_j$. The local isomorphism $\vphi_{j,\infb}$ yields a
$G$-equivariant isomorphism
\[
    (\E_j\times_{\D_j}\B_j)\otimes_R k \;\cong\;
      \Ind_{G_j}^G(\Yb_{j,\infb}) \;\lpfeil{\vphib_{j,\infb}}\;
      \Ind_{G\dagg,\lambda}^G(\Yb_{j,\infb}\dagg) \;\cong\;
        (\V\times_{\U}\B_j)\otimes_R k
\]
of \'etale $\Xb_{j,\infb}\dagg$-schemes. It lifts canonically to a
$G$-equivariant isomorphism
\[
     \vphi_j: \E_j\times_{\D_j}\B_j \;\liso\; \V\times_{\U}\B_j
\]
of \'etale formal $\B_j$-schemes.  By a theorem of Ferrand--Raynaud
\cite{FerrandRaynaud} (see also \cite{RachelLuminy}), there exists a
formal $\hat{X}_R$-scheme $\Y$ with $G$-action such that
$\E_j=\Y\times_{\hat{X}_R}\D_j$ and $\V=\Y\times_{\hat{X}_R}\U$. By
Grothendieck's Existence Theorem, $\Y$ is the formal completion of a
projective $R$-curve $Y_R$. By construction, the natural map
$f_R:Y_R\to X_R$ is a lift of $\fb$ which is mapped by
\eqref{patcheq2} to the patching datum $(\lambda;\vphib_{j,\infb})$.
This completes the proof of the theorem.
\end{proof}

\subsection{Some applications of Theorem \ref{liftthm}} 
                \label{corollaries}

Let $k$ be an algebraically closed field of characteristic $p>0$. We
denote by $K_0$ the fraction field of $W(k)$ and by $\Kb$ an algebraic
closure of $K_0$. Theorem \ref{liftthm} together with Proposition
\ref{patchprop2} implies:

\begin{cor} \label{liftcor0}
  Let $\fb:\Yb\to\Xb$ be a special $G$-map defined over $k$. Then
  there exists a three point $G$-cover $f:Y\to X:=\PP^1_{\Kb}$ whose
  stable reduction is isomorphic to $\fb$.
\end{cor}

This corollary is used in \cite{stabredmod} to determine the stable
reduction of all three point covers with Galois group ${\rm
PSL}_2(p)$. In particular, this gives a new proof for the result of
Deligne--Rapoport on the reduction of the modular curves $X_0(p)$ and
$X_1(p)$. 

In addition to the existence of three point covers with given stable
reduction, Theorem \ref{liftthm} and Proposition \ref{patchprop2}
yield interesting arithmetic information about a given three point
cover. 

\begin{cor} \label{liftcor1}
  Let $f:Y\to X=\PP^1_{\Kb}$ be a three point Galois cover with mild
  reduction. Let $(m_j,h_j)_{j\in B}$ be the signature of the stable
  reduction of $f$. The cover $f$ has stable reduction over the unique
  tame extension $K\stab/K_0$ of degree
  \[
        N \;:=\; (p-1)\cdot\mathop{\rm lcm}_{j\not\in B\wild}(h_j).
  \]
  In particular, any three point Galois cover $f:Y\to X$ whose degree
  is not divisible by $p^2$ can be defined over a tame extension of
  $K_0$.
\end{cor}

\begin{proof} Let $f_R:Y_R\to X_R$ be a stable model of $f$ and
$\fb:\Yb\to\Xb$ the stable reduction. If $\fb$ is not exceptional,
then the corollary follows from Proposition \ref{auxmonprop} and
Theorem \ref{liftthm}.  We may therefore assume that $\fb$ is
exceptional. We will also assume that we are in Case (2) of
Definition \ref{exceptdef}, leaving Case (1) for the reader. By
Proposition \ref{exceptprop}, the curve $\Xb$ consists of the original
component and a single new tail $\Xb_{j_0}$.

Let $f_R\aux:Y_R\aux\to X_R$ denote the auxiliary cover associated to
$f_R$, see \S \ref{patch1}. The Galois group of $f_R\aux$ is a
subgroup $G_0\subset G$, which is a {\em direct} product $G_0=I_0\times
H_0$. Here $I_0\cong\ZZ/p$ is the inertia group of a component $\Yb_0$
of $\Yb$ over $\Xb_0$ and $H_0$ is the Galois group of the cover
$\Zb_0=\Yb_0^{(p)}\to\Xb_0$. Let $W_R:=Y_R\aux/H_0$. Then $W_R\to X_R$
is a stable model of the $p$-cyclic Galois cover $W_K\to X_K=\PP^1_K$,
which is ramified at $0$, $1$ and $\infty$. Therefore, the curve $W_K$
has an equation of the form
\[
     y^p \;=\; f(x):=c\,x^a(x-1)^b,
\]
with integers $a,b$ such that $a,b,a+b\not\equiv 0\mod{p}$ and a unit
$c\in R^\times$. After enlarging the field $K$ if necessary and a
change of coordinates we may assume that $c\in W(k)^\times$ and
$f(a/(a+b))\in (W(k)^\times)^p$. It is known from
\cite{ColemanMcCallum} (see also \cite{Lehr01}) that the model $X_R$
which leads to the stable reduction of $W_K\to X_K$ is the blowup of
$\PP^1_R$ in the disk
\[
      D \;:=\; \{\,x \mid \val(x-\frac{a}{a+b}) \;\geq\; 
             \frac{p\val{p}}{2(p-1)} \;\}.
\]
Moreover, the stable model $W_R\to X_R$ descents to a stable model
over the ring
\[
     R' \;:=\; W(k)[\,p^{1/2(p-1)},\,f(\frac{a}{a+b})^{1/p}\,] \;=\;
               W(k)[\,p^{1/2(p-1)}\,].
\]
Note that the (admissible) $H_0$-cover $\Yb\aux\to\Wb:=W_R\otimes_Rk$
is ramified at the unique singular point of $\Wb$ of order
$m=h/2$. Therefore, the stable model $f_R\aux:Y_R\aux\to X_R$ descents
to a stable model over the ring $R'':=W(k)[p^{1/h(p-1)}]$. Using the
method of the proof of Theorem \ref{liftthm}, it is not hard to show
that $f_R:Y_R\to X_R$ descents to a stable model over $R''$. This
finishes the proof of the corollary.
\end{proof}

The last corollary can be refined in order to obtain a formula for the
ramification index of the field of moduli of a three point cover in
terms of its stable reduction.

\begin{cor} \label{liftcor2}
  Let $\fb:\Yb\to\Xb$ be a special $G$-map, with signature
  $(m_j,h_j)_{j\in B}$. Let $P\subset G$ denote a $p$-Sylow of $G$ and
  let $n:=[N_G(P):C_G(P)]$ be the index of the centralizer of $P$
  inside the normalizer of $P$. Recall that $\tilde{L}(\fb)$ denotes
  the set of isomorphism classes of three point $G$-covers defined
  over $\Kb$ whose stable reduction is isomorphic to $\fb$. There
  exists a positive integer $n'$ dividing $n$ which depends only on
  $\fb$ such that
  \[
      |\tilde{L}(\fb)| \;=\; \frac{p-1}{n'}
        \prod_{j\not\in B\wild}\frac{h_j}{|\Aut_G(\fb_j)|}.
  \]
  Let $f:Y\to X$ be a lift of $\fb$. Then
  the field of moduli $K\inn$ of $f$ (relative to the extension
  $\Kb/K_0$) is the unique tame extension of $K_0$ of degree
  \[
      N' \;:=\; \frac{p-1}{n'}
        \mathop{\rm lcm}_{j\not\in B\wild}\frac{h_j}{|\Aut_G(\fb_j)|}.
  \]
\end{cor}  
 
\begin{proof} We have seen in \S \ref{lifts2} that
$\tilde{L}(\fb)=L(\fb)/\Aut_G^0(\fb)$. Moreover, the stabilizer in
$\Aut_G^0(\fb)$ of any element of $L(\fb)$ is equal to $C_G$. Let $n'$
be the order of the image of the homomorphism $\Aut_G^0(\fb)\to G/C_G$
which is defined by Lemma \ref{autflem}. It is easy to see that $n'$
divides $n$. We see that Corollary \ref{liftcor2} follows from
Proposition \ref{auxmonprop}, Theorem \ref{liftthm} and Lemma
\ref{autflem}.
\end{proof}

\subsection{Genus zero dessins of degree $p$}

Let $p$ be a prime and $G\subset S_p$ be a primitive transitive
permutation group of degree $p$. We write $N\subset G$ for the
stabilizer of $1$. Let $f:Y\to X:=\PP^1_{\bar{\QQ}}$ be a three point
Galois cover, with Galois group $G$. Then $g:Z:=Y/N\to X$ is a three
point cover of degree $p$, and $f:Y\to X$ is the Galois closure of
$g$. We make the following assumptions.
 
\begin{ass} \label{dessinsass}
\begin{enumerate}
\item
  The subgroup $N\subset G$ is self-normalizing and $G$ is
  center-free.
\item
  The curve $Z$ has genus $0$.
\end{enumerate}
\end{ass}

The first assumption is made only for simplicity. It ensures that the
$G$-cover $f:Y\to X$ and the (non-Galois) cover $g:Z\to X$ have the same field
of moduli, and that this field of moduli is also the minimal field of
definition for both covers, see e.g\ \cite{DebDou1}. The second assumption is
more restrictive and implies that the stable reduction of $f:Y\to X$ to
characteristic $p$ is rather easy to describe.

Let $K$ be the field of moduli of $f:Y\to X$, relative to the
extension $\bar{\QQ}/\QQ$. Let $\p$ be a prime ideal of $K$ dividing
$p$. There exists a finite extension $K'/K$ and a prime $\p'$ of $K'$
over $\p$ such that $f$ has a $K'$-model $f_{K'}:Y_{K'}\to X_{K'}$
which extends to a stable model over the local ring $\OO_{K',\p'}$. Let
$\fb:\Yb\to\Xb$ be the stable reduction of $f_K'$ at $\p'$. One checks
that $\fb:\Yb\to\Xb$ depends only on the field $K$ and the prime $\p$.
Hence we may call $\fb:\Yb\to\Xb$ the stable reduction of $f:Y\to X$
at $\p$.

Let us assume that $f:Y\to X$ has bad reduction at $\p$. Note that $p$
strictly divides the order of $G$. Therefore, by Theorem \ref{thm1} the stable
reduction $\fb:\Yb\to\Xb$ is a special $G$-map, associated to a special
$G$-deformation datum $(\Zb_0,\omega_0;\Yb_j)$, see \S \ref{Gdefdat3}.  In the
following, we use freely the notation introduced in \S \ref{Gdefdat3}.

\begin{lem} \label{dessinslem}
\begin{enumerate}
\item
  The character $\chi:H_0\to\FF_p^\times$ is injective. In particular,
  $H_0\cong\ZZ/m$ with $m|(p-1)$.
\item
  There are no new tails.
\item
  For $j\in B\prim$, let $x_j\in\{0,1,\infty\}$ denote the branch
  point of $f:Y\to X$ which specializes to the primitive tail
  $\Xb_j$. The ramification invariant of the primitive tail
  cover $\fb_j:\Yb_j\to\Xb_j$ is given by the formula
  \[
       \sigma_j \;=\; \frac{|g^{-1}(x_j)|-1}{p-1}.
  \]
\end{enumerate}
\end{lem}

\begin{proof}
Part (1) is an immediate consequence of the fact that the $p$-cyclic
subgroup $I_0$ of $S_p$ is self-centralizing. 

To prove (2), suppose that $f_j:\Yb_j\to\Xb_j$ is a new tail, with
Galois group $G_j$. Let $N_j:=N\cap G_j$ and $\Zb_j:=\Yb_j/N_j$. It
follows from Assumption \ref{dessinsass} (2) that the curve $\Zb_j$
has genus $0$. Since $N$ has index $p$ in $G$ and $G_j$ contains a
$p$-Sylow of $G$, the index of $N_j$ in $G_j$ is also $p$. We see that
$\gb_j:\Zb_j\to\Xb_j$ is a cover of degree $p$ between curves of genus
$0$, totally branched at $\infb_j\in\Xb_j$ and \'etale elsewhere. This
means that in suitable coordinates $\gb_j$ is given by a monic
polynomial $q\in k[z]$ of degree $p$ such that the derivative $q'$ has
degree $0$. Hence $q=z^p+az+b$, and after a further normalization we
get $q=z^p-z$. It follows that the tail cover $\fb_j:\Yb_j\to\Xb_j$ is
an Artin--Schreier cover with conductor $1$, i.e.\ $\sigma_j=1$. This
contradicts the results of \S \ref{threepoints}. We conclude that
$\fb:\Yb\to\Xb$ has no new tails, proving (2).

For $j\in B\prim$, let $\gb_j:\Zb_j\to\Xb_j$ be as before. The same argument
as before shows that $\gb_j$ is a cover between smooth curves of genus $0$ of
degree $p$, totally ramified at $\infb_j\in\Xb_j$. But now $\gb_j$ is tamely
ramified at $\xb_j$ and \'etale over $\Xb_j-\{\infb_j,\xb_j\}$ (recall that
$\xb_j\in\Xb_j$ denotes the specialization of the branch point
$x_j\in\{0,1,\infty\}$). Since the stable model of $f$ is a tame cover in a
neighborhood of $\Xb_j-\{\infb_j\}$ we have
$|\gb_j^{-1}(\xb_j)|=|g^{-1}(x_j)|$. Now (3) follows from a calculation with
the Riemann--Hurwitz formula and the fact that $g(\Zb_j)=0$.  
See e.g.\ \cite{RRR}, Lemma 2.3.1.
\end{proof}

The lemma shows that the signature of $\fb:\Yb\to\Xb$ is a tripel
$(\sigma_1,\sigma_2,\sigma_3)$ of rational numbers, with $\sigma_j=h_j/m_j$,
$0\leq h_j<m_j$ and $(h_j,m_j)=1$. It satisfies the vanishing cycle formula
\begin{equation} \label{dessineq1}
  \sigma_1 + \sigma_2 + \sigma_3 \;=\; 1.
\end{equation}
The special deformation datum $(\Zb_0,\omega_0)$ is essentially determined by
the signature. Indeed, the curve $\Zb_0$ is isomorphic to the smooth
projective model of the plane curve with equation
\begin{equation} \label{dessinseq2}
       z^m \;=\; x^{h_1}(x-1)^{h_2},
\end{equation}
where $m$ is the least common multiple of the $m_j$, and the differential
$\omega_0$ has the form
\begin{equation} \label{dessinseq3}
   \omega_0 \;=\; \epsilon\,\frac{z\,\diff x}{x(x-1)},
\end{equation}
for some $\epsilon\in\FF_p^\times$. Conversely, given a tripel
$(\sigma_1,\sigma_2,\sigma_3)$ satisfying all of the above conditions, one
checks that \eqref{dessinseq2} and \eqref{dessinseq3} define a special
deformation datum $(\Zb_0,\omega_0)$. Therefore, the existence of a special
deformation datum with a given signature imposes no restriction on the stable
reduction $\fb:\Yb\to\Xb$.


\begin{exa} \label{dessinsexa1}
  Let $p=7$ and $G:=S_7$.  Suppose furthermore that the branch cycle
  description of the cover $g:Z\to X$ is $(6,6,\text{$2$-$2$})$ (each entry
  designates a conjugacy class of $S_7$, e.g.\ `$\text{$2$-$2$}$' denotes the
  product of two $2$-cycles). From the table computed by Malle \cite{Malle94}
  we see that there are exactly $4$ non-isomophic three point covers with
  monodromy group $S_7$ and this branch cycle description. Moreover, the
  absolute Galois group $\Gal(\bar{\QQ}/\QQ)$ permutes these four covers
  transitively, with permutation group $S_4$. Therefore, the field of moduli
  of the $G$-cover $f:Y\to X$ is a degree $4$ extension $K/\QQ$.
  
  Let $\p$ be a prime ideal of $K$ over $7$. Suppose that $f:Y\to X$
  has bad reduction at $\p$ and let $\fb:\Yb\to\Xb$ denote the stable
  reduction. It follows from Lemma \ref{dessinslem} (3) that the
  signature of $\fb$ is equal to $(1/6,1/6,2/3)$. For $j=1,2$, the
  Riemann--Hurwitz formula implies that the Galois group $G_j$ of the
  primitive tail cover $\fb_j:\Yb_j\to\Xb_j$ is isomorphic to the
  semi-direct product $7:6$. It follows that
  $\Yb_j\to\Xb_j$ is given by the two equations
  \begin{equation} \label{dessinseq4}
        z^6 \;=\; x, \qquad y^7-y \;=\; z.
  \end{equation}
  The third tail cover $\fb_3:\Yb_3\to\Xb_3$ has inertia invariant
  $\sigma_3=2/3$ and branch cycle description $\text{$2$-$2$}$ at the tame
  branch point.  We say that $\fb_3$ is a {\em primitive tail cover of type
    $(2/3,\text{$2$-$2$})$}.  The Galois group $G_3$ of $\fb_3$ is a subgroup
  of $A_7$ containing a semidirect product $7:3$ and a permutation of type
  $\text{$2$-$2$}$. By elementary group theory, $G_3$ is isomorphic to $A_7$
  or to $L_2(7):={\rm PSL}_2(7)$.
  
  It is a priori not clear how many non-isomorphic primitive tail covers of
  type $(2/3,\text{$2$-$2$})$ exist, and what their Galois groups are. We
  claim that there exists in fact a {\em unique} primitive tail cover of type
  $(2/3,\text{$2$-$2$})$, with Galois group $L_2(7)$. The existence of one
  such cover is proved in \cite{stabredmod}. (It is obtained by reducing a
  certain three point cover with Galois group $L_2(7)$.) To prove uniqueness,
  we use our lifting result, Theorem \ref{liftthm}.
  
  Given $\fb_3:\Yb_3\to\Xb_3$ as above, one shows that there exists, up to
  isomorphism, a unique special $G$-map $\fb:\Yb\to\Xb$ with signature
  $(1/6,1/6,2/3)$ whose third tail cover is isomorphic to
  $\fb_3:\Yb_3\to\Xb_3$. By Corollary \ref{liftcor2} there exist exactly
\[
         N' \;:=\; \frac{12}{n'\cdot|\Aut_{G_3}(\fb_3)|}
\]
nonisomorphic three point covers $f':Y'\to X'$ whose stable reduction (at a
given prime $\tilde{\p}$ of $\bar{\QQ}$) is isomorphic to $\fb$.  Moreover,
all of these covers are conjugate under the action of the inertia group of
$\tilde{\p}$.  Using Lemma \ref{autflem}, Remark \ref{autfrem} and
Equation \eqref{dessinseq4} one shows that $n'=3$. Thus, $N'$ is equal to
either $2$ or $4$.

As we have mentioned before, all three point covers with Galois group $S_7$
and branch cycle description $(6,6,\text{$2$-$2$})$ are conjugate. Hence we
may assume that $f=f'$. We conclude that for each primitive tail cover $\fb_3$
of type $(2/3,\text{$2$-$2$})$ there exists a prime $\p$ in $K$ above $7$ with
$e(\p/7)\in\{2,4\}$. This prime $\p$ is uniquely determined by the property
that the stable reduction of $f:Y\to X$ at $\p$ has a primitive tail cover
isomorphic to $\fb_3$. By Malle's table, we have $7=\p_1^2\p_2\p_3$ in $K$.
Therefore, our analysis of the stable reduction implies:
\begin{itemize}
\item The cover $f:Y\to X$ has bad reduction at $\p_1$ and good reduction at
  $\p_2$ and $\p_3$. The stable reduction of $f:Y\to X$ at $\p_1$ occurs after
  an extension which is ramified at $\p_1$ of order $2$.
\item There exists a unique primitive tail cover of type
  $(2/3,\text{$2$-$2$})$. Its Galois group is $L_2(7)$.
\end{itemize}
We remark that good reduction at $\p_2$ and $\p_3$ does not follow from the
results of \cite{Raynaud98}, because $S_7$ has only one conjugacy class of
elements of order $7$. 
\end{exa}

\begin{exa} \label{dessinsexa2}
  As before, let $p=7$ and $G=S_7$. But now suppose that the branch cycle
  description of the cover $f:Y\to X$ is $(\text{$2$-$3$},\text{$2$-$3$},7)$.
  By Malle's table there are 9 such covers, up to isomorphism. Note that they
  all have bad reduction to characteristic $7$, because of the ramification
  index $7$.
  
  It follows from Lemma \ref{dessinslem} (3) that the signature of
  $\fb:\Yb\to\Xb$ is equal to $(1/2,1/2,0)$. For $j=1,2$, the Galois
  group $G_j$ of the tail cover $\fb_j:\Yb_j\to\Xb_j$ contains a
  $7$-cycle and a $2$-cycle. Hence, by a theorem of Jordan we have
  $G_j=S_7$.
  
  We claim that there exists, up to isomorphism, a unique primitive tail cover
  $\fb_1:\Yb_1\to\Xb_1$ of type $(1/2,\text{$2$-$3$})$ with Galois group
  $S_7$. We have already shown its existence. To prove uniqueness, we use
  again Theorem \ref{liftthm}.
  
  Let us choose a point $\eta_1\in\Yb_1$ above the wild branch point
  $\infb_1\in\Xb_1$. Since $\sigma_1=1/2$ the inertia group at $\eta_1$ is a
  dihedral group $G_0=I_0\rtimes H_0$, with $|I_0|=7$ and $H_0=|2|$.
  Recall that the triple $(\fb_1,\eta_1,H_0)$ is called a pointed $S_7$-tail
  cover. Let $\fb_2:\Yb_2\to\Xb_2$ be a copy of the $S_7$-tail cover $\fb_1$.
  It is easy to see that there are exactly three points on $\Yb_2$ above
  $\infb_2$ with inertia group $G_0$. Let $\eta_2$ be one of them. Then
  $(\fb_2,\eta_2,H_0)$ is another pointed $S_7$-tail cover. It is isomorphic
  to $(\fb_1,\eta_1,H_0)$ (by an isomorphism commuting with the action of
  $S_7$) if and only if $\eta_1=\eta_2$. Finally, let $(\Zb_0,\omega_0)$ be
  the (unique) special deformation datum with signature $(1/2,1/2,0)$. One
  checks that $(\Zb_0,\omega_0,\Yb_1,\Yb_2)$ gives rise to a special
  $G$-deformation datum and hence to a special $G$-map $\fb:\Yb\to\Xb$, see \S
  \ref{Gdefdat3}.  Moreover, by varying the point $\eta_2$ we obtain exactly
  three non-isomorphic special $G$-maps with signature $(1/2,1/2,0)$ and whose
  primitive tail covers are both isomorphic to $\fb_1:\Yb_1\to\Xb_1$.
  
  Choose one of the three special $G$-maps $\fb:\Yb\to\Xb$ constructed from
  $\fb_1:\Yb_1\to\Xb_1$.  By Corollary \ref{liftcor2} there exist exactly
  $N'=6/n'$ non-isomorphic three point covers $f':Y'\to X$ with Galois group
  $S_7$ whose stable reduction (at a given place $\tilde{\p}$ above $7$) is
  isomorphic to $\fb$. Moreover, they form a single orbit under the action of
  the inertia group at $\tilde{\p}$. Using Remark \ref{autfrem} we see that
  $2|n'$. Using the fact that $S_7$ has no outer automorphism one further
  shows that $n'=2$ and hence $N'=3$.  Therefore, for each primitive tail
  cover $\fb_1:\Yb_1\to\Xb_1$ we obtain, by our lifting result, exactly $9$
  non-isomorphic three point covers with Galois group $S_7$ and branch cycle
  description $(\text{$2$-$3$},\text{$2$-$3$},7)$. However, there are only $9$
  such covers altogether.  This proves that $\fb_1:\Yb_1\to\Xb_1$ is indeed
  unique, up to isomorphism. As a by-product, we have also shown that the
  absolute ramification index of any prime ideal $\p$ dividing $7$ in the
  field of moduli of $f:Y\to X$ is equal to $3$.
  
  By \cite{Malle94} the set of $9$ isomorphism classes form two Galois
  orbits, of length $3$ and $6$ respectively, corresponding to extensions
  $K_1/\QQ$ and $K_2/\QQ$. In $K_1$ we have $7=\p_1^3$ and in $K_2$ we have
  $7=\p_2^3\p_3^3$. This confirms our prediction. 
\end{exa}

\begin{rem}
  One can use Malle's table \cite{Malle94} to verify our results in many more
  examples. In particular, one gets a numerical verification of Theorem
  \ref{introthm1} of the Introduction for all genus zero dessins of degree
  $\leq 10$. One also sees that the condition that $p$ strictly divides the
  order of $G$ is necessary for Theorem \ref{introthm1} to hold.  For
  instance, one finds several three point covers of degree $n\geq 10$, with
  monodromy group $S_n$ or $A_n$, such that $5$ is wildly ramified in the
  field of moduli.
\end{rem}  

\begin{rem}
The relation between the existence of three point covers with bad reduction on
the one hand and tail covers with low ramification invariant on the other hand
is exploited more systematically in \cite{RRR}. In \cite{Zapponi02} Zapponi
studies the stable reduction of genus zero dessins of degree $p$ which are
totally ramified at one point, by analyzing the defining equation. His results
confirm our results on dessins of genus $0$. 
\end{rem}

\providecommand{\bysame}{\leavevmode\hbox to3em{\hrulefill}\thinspace}


\end{document}

%% file: stablepic.tex

\begin{picture}(30,16)

\put(6,4){\line(1,0){16}}
\put(6,10){\line(1,0){16}}
\put(6,14){\line(1,0){16}}

\put(10,3){\line(0,1){4}}
\put(12,3){\line(0,1){4}}
\put(16,3){\line(0,1){4}}
\put(19,3){\line(0,1){4}}

\put(10,9){\line(0,1){6}}
\put(12,9){\line(0,1){6}}
\put(16,9){\line(0,1){6}}
\put(19,9){\line(0,1){6}}

\put(8,4){\circle*{0.3}}
\put(10,6){\circle*{0.3}}
\put(12,6){\circle*{0.3}}

\put(8,2){\makebox(0,0){$\scriptstyle 0$}}
\put(10,2){\makebox(0,0){$\scriptstyle 1$}}
\put(12,2){\makebox(0,0){$\scriptstyle \infty$}}

\put(17,5){\makebox(1,1){$\cdots$}}
\put(17,12){\makebox(1,0){$\cdots$}}

\put(7,11.9){\makebox(0,1){$\vdots$}}
\put(14,11.9){\makebox(0,1){$\vdots$}}
\put(21,11.9){\makebox(0,1){$\vdots$}}

\put(1.5,11.9){\makebox(1,1){$\Yb$}}
\put(1.5,3.5){\makebox(1,1){$\Xb$}}
\put(2,10){\vector(0,-1){3.5}}



\put(24.2,4){\makebox(0,0){$\Xb_0=\PP^1_k$}}

\end{picture}